\documentclass[11pt,reqno]{amsart}
\usepackage{amsmath,amsthm,amssymb,amsfonts,mathtools,bm}
\usepackage{mathrsfs}
\usepackage{fullpage}
\usepackage{cite}
\usepackage[colorlinks=true,citecolor=blue,linkcolor=blue,urlcolor=blue]{hyperref}
\usepackage[colorlinks=true,citecolor=blue,linkcolor=blue,urlcolor=blue]{hyperref}
\pdfmapfile{+cm.map}
\pdfmapfile{+cmextra.map}
\pdfmapfile{+symbols.map}

\newcommand{\bea}{\begin{eqnarray}}
\newcommand{\eea}{\end{eqnarray}}
\newcommand{\bna}{\begin{eqnarray*}}
\newcommand{\ena}{\end{eqnarray*}}

\numberwithin{equation}{section}

\theoremstyle{plain}
\newtheorem{lem}{Lemma}[section]
\newtheorem{theorem}[lem]{Theorem}

\newtheorem{proposition}[lem]{Proposition}

\theoremstyle{definition}

\newtheorem{remark}{Remark}

\newcommand*{\dif}{\mathop{}\!\mathrm{d}}

\begin{document}

\title{The Variance of Isobaric Automorphic Coefficients in Arithmetic Progressions}
\author{Qingfeng Sun}
	\address{School of Mathematics and Statistics, Shandong University, Weihai\\Weihai, Shandong 264209, China}
	\address{State Key Laboratory of Cryptography and Digital Economy
Security \\ Shandong University, Weihai, 264209, China}
	\email{qfsun@sdu.edu.cn}

    \author{Qizhi Zhang}
	\address{School of Mathematics and Statistics, Shandong University, Weihai\\Weihai, Shandong 264209, China}
    \email{qzzhang@mail.sdu.edu.cn}

\thanks{Q. Sun was partially supported by the National Natural Science Foundation of China
(Grant No. 12471005) and
the Natural Science Foundation of Shandong Province (Grant No. ZR2023MA003). Q. Zhang was supported by the Program of China Scholarship Council (Grant No. 202606220052).}

\date{}

\begin{abstract}
Let \(\pi\) be a fixed isobaric representation of \(\mathrm{GL}_m(\mathbb A_{\mathbb Q})\)
whose \(r\) unitary cuspidal constituents have degree at least two and are
pairwise inequivalent up to twists by \(|\det|^{it}\), \(t\in\mathbb R\).
For normalized smooth weights \(\Phi,\Psi\), the variance of the coefficients
\(\lambda_\pi(n)\) in reduced residue classes satisfies
\[
\begin{aligned}
 V_\pi(Q,X)
 &= \sum_{q\ge1}\Phi(q/Q)
    \sum_{\substack{a\bmod q\\(a,q)=1}}
    \left|\sum_{n\equiv a\bmod q}\lambda_\pi(n)\Psi(n/X)
      -\frac{1}{\phi(q)}\sum_{(n,q)=1}\lambda_\pi(n)\Psi(n/X)\right|^2 \\
 &= QX P_{r-1}(\log X)+o(QX),
\end{aligned}
\]
where \(P_{r-1}\) is an explicit polynomial of degree \(r-1\).
The formula holds uniformly as \(Q\to\infty\) for
\(Q^\varepsilon\le X\le Q^{1+2/(m+2)-\varepsilon}\). Under the generalized
Ramanujan conjecture and the generalized Riemann hypothesis, the asymptotic
formula holds for \(Q^\varepsilon\le X\le Q^{2-\varepsilon}\).
\end{abstract}

\maketitle

\section{Introduction}

The quantitative study of arithmetic functions in residue classes finds its
origin in the distribution of prime numbers.
Put
\[
 \psi(q,a;X)=\sum_{\substack{n\le X\\n\equiv a(\!\!\bmod q)}}\Lambda(n).
\]
The Barban--Davenport--Halberstam theorem \cite{Barban,DH}, in the strengthened form
due to Gallagher \cite{Gallagher67}
(see also \cite[Theorem~2]{Hooley75}), states that,
for every fixed $A>0$,
\[
\sum_{q\le Q}
\sum_{\substack{a(\!\bmod q)\\(a,q)=1}}
\left|\psi(q,a;X)-\frac{X}{\phi(q)}\right|^2
\ll_A XQ\log X+\frac{X^2}{(\log X)^A},
\]
uniformly for $X\ge3$ and $1\le Q\le X$.
This variance estimate illustrates the strength of averaging
over moduli and reduced residue classes. In applications that permit such averaging, estimates of this kind
typically yield substantially sharper aggregate bounds
than those obtained by simply summing the standard bounds under the generalized Riemann hypothesis (GRH) for individual
arithmetic progressions; see, for example,
Zhang \cite{Zhang14}
and Heath-Brown and Li \cite{HBL17}.

Beyond the distribution of primes, variance problems in
arithmetic progressions have also been extensively studied
for other arithmetic sequences, notably the divisor functions.
For a fixed integer $k\ge2$, let $d_k(n)$ denote the
$k$-fold divisor function, defined by
\[
d_k(n):=\#\{(a_1,\ldots,a_k)\in\mathbb{Z}_{\ge1}^k:
a_1\cdots a_k=n\}.
\]
The distribution of $d_k(n)$ in arithmetic progressions was studied
in the early work of Lavrik \cite{Lavrik1965,Lavrik1966}.
Motohashi \cite{Motohashi1973} subsequently established an asymptotic
formula for the variance of $d_2(n)$, averaged over moduli.
For further variance bounds and related mean-square
asymptotics for divisor functions and their generalizations, see
\cite{Blomer2008,HarperSoundararajan2017,LauZhao2012,FGKM2014,%
PongsriiamVaughan2018,deLaBretecheFiorilli2020,KowalskiRicotta2014,%
PongsriiamVaughan2021,Mastrostefano2022,%
Nguyen2024Individual,%
Parry2025} and the references therein.
Of these developments, the work of Rodgers and Soundararajan
\cite{RodgersSoundararajan2018} is particularly noteworthy.
Their work was motivated in part by the
function-field results of Keating, Rodgers, Roditty-Gershon
and Rudnick \cite{KRRR}, which relate analogous variances over
$\mathbb{F}_q[T]$ to unitary matrix integrals and suggest
a piecewise-polynomial profile for the normalized variance
over the integers.
For fixed nonnegative weights
$\Phi,\Psi\in C_c^\infty((0,\infty))$, normalized by
$\int_0^\infty\Phi(u)\,du=\int_0^\infty\Psi(u)^2\,du=1$,
they considered
\[
 V_k(Q,X)
 =
 \sum_{q\geq1}\Phi\!\left(\frac qQ\right)
 \sum_{\substack{a\bmod q\\(a,q)=1}}
 \bigg|
  \sum_{n\equiv a\!\!\!\!\pmod q}
  d_k(n)\Psi\left(\frac nX\right)
  -
  \frac1{\phi(q)}
  \sum_{(n,q)=1}
  d_k(n)\Psi\!\left(\frac nX\right)
 \bigg|^2.
\]
Using the asymptotic large sieve, they proved
the conjectured asymptotic for $V_k(Q,X)$ uniformly,
for each fixed $k\ge2$ and $\varepsilon>0$, in the range
\[
 Q^\varepsilon\leq X\leq Q^{1+2/k-\varepsilon}.
\]
Under the generalized Riemann hypothesis (GRH),
the same asymptotic holds uniformly for
$Q^\varepsilon\leq X\leq Q^{2-\varepsilon}$.

The purpose of the present paper is to study an automorphic analogue of
this variance problem.
Note that
\begin{equation*}
\sum_{n=1}^{\infty}\frac{d_k(n)}{n^s}
=\zeta(s)^k
=L(s,\mathbf{1}^{\boxplus k}),
\qquad \operatorname{Re}(s)>1,
\end{equation*}
where $\mathbf{1}^{\boxplus k}$ denotes the isobaric sum of $k$
copies of the trivial automorphic character of
$\mathrm{GL}_1(\mathbb A_{\mathbb Q})$.
This motivates us to replace the trivial constituents by
more general cuspidal automorphic representations
and to study analogous variance problems for the resulting
coefficients.

More precisely, let $\pi_j$, $1\le j\le r$, be fixed nontrivial
unitary cuspidal automorphic representations of
$\mathrm{GL}_{m_j}(\mathbb A_{\mathbb Q})$ with $m_j\geq 2$, and
assume throughout that $\pi_h\not\simeq \pi_j\otimes|\det|^{it}$
for all $h\ne j\text{ and }t\in\mathbb R.$
Set
\[
\pi=\pi_1\boxplus\cdots\boxplus\pi_r, \qquad m=\sum_{j=1}^r m_j.
\]
The isobaric $L$-function associated to $\pi$ is of degree $m$ and defined by
\[
L(s,\pi)
=\prod_{j=1}^{r}L(s,\pi_j)
=\sum_{n=1}^{\infty}\frac{\lambda_\pi(n)}{n^s},
\qquad \operatorname{Re}(s)>1,
\]
where $L(s,\pi_j)$ is the standard $L$-function associated to $\pi_j$.
Writing \(\lambda_{\pi_j}\) for the Dirichlet coefficients of
\(L(s,\pi_j)\), then $\lambda_\pi=
\lambda_{\pi_1}\ast\cdots\ast\lambda_{\pi_r}.$

Define the variance
\begin{align}\label{DeltaPiQX}
V_\pi(Q,X)
 :=
 \sum_{q\geq1}\Phi\!\left(\frac qQ\right)
 \sum_{\substack{a\!\!\!\!\! \!\pmod q\\(a,q)=1}}\;
 \Bigl| \sum_{n\equiv a\!\!\!\!\! \!\pmod q}
 \lambda_\pi(n)\Psi\!\left(\frac nX\right)-
 \frac{1}{\phi(q)}\sum_{(n,q)=1}\lambda_\pi(n)\Psi\left(\frac{n}{X}\right)\Bigr|^2,
\end{align}
where \(\Phi,\Psi\in C_c^\infty((0,\infty);\mathbb R)\) satisfy
\[
 \Phi\ge 0,\qquad
 \int_0^\infty \Phi(x)\dif x
 =
 \int_0^\infty \Psi(x)^2\dif x
 = 1.
\]

Our main result is the following.
\begin{theorem}\label{mainresult1}
For every $\varepsilon>0$ and
\[
 Q^\varepsilon\leq X\leq Q^{1+2/(m+2)-\varepsilon},
\]
we have
\[
V_\pi(Q,X)=QX P_{r-1}(\log X)+o(QX),
\]
where $P_{r-1}$ is a polynomial of degree $r-1$.  Its leading coefficient is
\begin{align}\label{eq:leading-constant}
 \frac{B_{\pi,0}(1)H_\pi(1)}{(r-1)!}
 \prod_{h=1}^r
 \operatorname*{Res}_{s=1}
 L(s,\pi_h\times\widetilde\pi_h)
 \prod_{\substack{h,j=1\\h\ne j}}^r
 L(1,\pi_h\times\widetilde\pi_j),
\end{align}
where $H_\pi(s)$ is defined in \eqref{eq:H-factor}, and
\begin{align}\label{eq:B-one-intro}
B_{\pi,0}(1)=\prod_p\bigg\{1+
\frac{1}{p}
 \bigg(\Big(1+\sum_{\ell\geq1}\frac{|\lambda_{\pi}(p^\ell)|^2}{p^\ell}\Big)^{-1}-1\bigg)
 \bigg\}.
\end{align}
The Euler product in \eqref{eq:B-one-intro} converges absolutely.
\end{theorem}

\begin{remark}
The proof yields more information than is recorded in the preceding
asymptotic formula. For simplicity, we state the explicit remainder below
only when \(m>2\). Set $m^*=\max_{1\le j\le r}m_j$
and
\[
\theta=\frac{1}{2}-\frac{1}{m^{*2}+1}<\frac{1}{2}.
\]
The proof gives
\[
V_\pi(Q,X)
 =
 QX P_{r-1}(\log X)+R(Q,X),
\]
where
\[
 P_{r-1}(Y)
 =
 \sum_{k=0}^{r-1}\kappa_{\pi,k}Y^k
\]
is a polynomial of degree \(r-1\), and
\[
\begin{aligned}
 R(Q,X)
 =O_{\pi,\Phi,\Psi,\varepsilon}\biggl(
 Q^{2\theta+\varepsilon}X^{1+\varepsilon}
 &+QX^{1-1/(m^{*2}+1)+\varepsilon}
 \\&+X^{5/3+\varepsilon}
   \left(\frac{X}Q\right)^{(m-2)/6}+\frac{X^{2+\varepsilon}}{Q}
 +X^{1+\varepsilon}
   \left(\frac{X}Q\right)^{m/2}
 \biggr).
\end{aligned}
\]
Consequently, on every fixed closed subrange of the range specified in the
preceding theorem,
\[
 R(Q,X)=o(QX)
\]
uniformly as \(Q\to\infty\).

The residue calculation in the proof determines every coefficient
\(\kappa_{\pi,k}\) explicitly, not merely the leading coefficient
\(\kappa_{\pi,r-1}\). Thus the terms
\[
 QX\sum_{k=0}^{r-2}\kappa_{\pi,k}(\log X)^k
\]
are explicit lower-order main terms rather than part of the genuine analytic
remainder \(R(Q,X)\). The sum is empty when $r=1$
\end{remark}

\begin{theorem}\label{mainresult2}
Assume the generalized Ramanujan conjecture for each $\pi_j$ and GRH for
$L(s,\pi_j\times\chi)$ and $L(s,\widetilde\pi_j\times\overline\chi)$
for every primitive Dirichlet character $\chi$, including the conductor-one character.
Then the conclusion of Theorem~\ref{mainresult1} holds uniformly for
$Q^\varepsilon\leq X\leq Q^{2-\varepsilon}$.
\end{theorem}

For the proof of Theorems \ref{mainresult1} and \ref{mainresult2}, the automorphic setting introduces several new difficulties.  The elementary
Euler products occurring in the divisor-function problem must be replaced by
factorizations involving Rankin--Selberg \(L\)-functions, while the
nonprincipal-character contribution requires hybrid mean-value estimates for
the twists \(L(s,\pi\otimes\chi)\).  The shape of the main term is governed by
the square-coefficient Dirichlet series
\[
 D_\pi(s)
 :=
 \sum_{n=1}^\infty\frac{|\lambda_\pi(n)|^2}{n^s}.
\]
We prove that
\begin{align*}
 D_\pi(s)
 &=
H_\pi(s)L(s,\pi\times\widetilde\pi),
\end{align*}
where \(H_\pi(s)\) is holomorphic near \(s=1\) and
\(H_\pi(1)\neq0\).  Since the \(\pi_j\) are pairwise inequivalent,
\(L(s,\pi_h\times\widetilde\pi_j)\) has a simple pole at \(s=1\) precisely
when \(h=j\).  It follows that \(D_\pi(s)\) has a pole of order \(r\) at
\(s=1\), giving rise to a main term that is polynomial of degree \(r-1\) in
\(\log X\).  By comparison,
\(\sum_n d_k(n)^2n^{-s}\) has a pole of order \(k^2\) at \(s=1\).

Our unconditional argument uses Jiang's estimate for Euler products of
higher-prime-power Rankin--Selberg coefficients
\cite[Theorem~1.3]{J25}. This provides the averaged control needed for  the
correction factors in Lemma~\ref{corollary:Jiang} without assuming the
generalized Ramanujan conjecture (GRC). Pointwise Satake bounds remain
necessary elsewhere, while the wider conditional range uses both GRH and GRC.

\subsection*{Structure of the paper} Section~2 records the preliminary results and estimates.
Section~3 reduces the main theorems to three estimates, proved in Sections~ 4, 5 and~7.
Section~6 supplies the Mellin-transform estimates and the factorization used in the off-diagonal
analysis; the latter is proved in Section~8.

\subsection*{AI Declaration Statement}
All of the ideas and all the calculations in this manuscript are due to the authors. This work is finished dating back to May-2026.

\section{Preliminaries}

We review here some standard facts about
automorphic $L$-functions. For much of these, we refer
the reader to
Brumley \cite[Section 1]{Brum}.

\subsection{$L$-functions}

Let $\mathcal{A}(m)$ be the set of
unitary cuspidal automorphic
representations of
$\mathrm{GL}_{m}(\mathbb{A}_{\mathbb{Q}})$ with unitary
central characters; here $m\ge2$.
For \(\mathrm{Re}(s)>1\),
the global $L$-function associated to
$\pi=\otimes_v\pi_{v}\in \mathcal{A}(m)$ is defined by
\bna
L(s,\pi)=\prod_{p}L(s,\pi_p)
=\sum_{n=1}^{\infty}\frac{\lambda_{\pi}(n)}{n^s},
\ena
where for each finite prime $p$,
the inverse of the local factor $L(s,\pi_{p})$ is a polynomial
in $p^{-s}$ of degree less than or equal to $m$, i.e.,
\bna\label{mu-def}
L(s,\pi_{p})^{-1}=\prod_{j=1}^{m} (1-\alpha_{\pi,j}(p)p^{-s}),
\ena
with $\{\alpha_{\pi,j}(p)\}_{j=1}^{m}$ the
Satake parameters of $\pi_p$. For \(\mathrm{Re}(s)>1\), define
\bna
L(s,\pi)^{-1}=\prod_p L(s,\pi_{p})^{-1}
:=\sum_{n=1}^{\infty}
\frac{\mu_{\pi}(n)}{n^s}.
\ena

Define the local factor at the infinite place
\[
 L(s,\pi_{\infty}):=\prod_{j=1}^{m}
 \Gamma_{\mathbb R}(s+\nu_{\pi,j}),
 \hspace{2em}
 \Gamma_{\mathbb R}(s)=\pi^{-s/2}\Gamma(s/2),
\]
where $\{\nu_{\pi,j}\}_{j=1}^{m}$ are
Langlands parameters of $\pi_{\infty}$.
The completed function
\[
 \Lambda(s,\pi)=q_{\pi}^{s/2}
 L(s,\pi_{\infty})L(s,\pi)
\]
can be analytically continued to the whole complex plane
and satisfies the functional equation
\[
 \Lambda(s,\pi)
 =\varepsilon(\pi)\Lambda(1-s,\widetilde\pi),
\]
where $q_{\pi}$ is the arithmetic conductor
of $\pi$, $\varepsilon(\pi)$ is a complex number
of modulus $1$,
and $\tilde\pi\in \mathcal{A}(m)$ is the contragredient of $\pi$.
At every unramified finite place,
the Satake parameters of $\widetilde\pi$ are the complex
conjugates of those
of $\pi$ and thus
$\lambda_{\widetilde\pi}(n)=\overline{\lambda_{\pi}(n)}$.

The bounds of Kim and Sarnak \cite{ks} for
$2 \leq m \leq 4$
and Luo et al. \cite{LRS1}, \cite{LRS2} for $m \geq 5$ are
\bea\label{theta}
|\alpha_{\pi,j}(p)|\le p^{\theta_{m}},\quad\qquad
\mathrm{Re}(\nu_{\pi,j})\ge -\theta_{m}
\eea
for all primes $p$ and $1\le j\le m$, where
\bea\label{theta-j}
\theta_2=\frac{7}{64},\quad \theta_3=\frac{5}{14},\quad
\theta_4=\frac{9}{22},\quad
\theta_{m}=\frac12-\frac1{m^2+1}\ (m\ge 5).
\eea
Moreover,
\bea\label{lambda-n}
 \bigl|\lambda_{\pi}(n)\bigr|
 \leq d_{m}(n)n^{\theta_{m}},
\eea
where $d_k(n)=
\sharp\{(a_1,a_2,\cdots,a_k)\in \mathbb{Z}_{>0}^k:a_1a_2\cdots a_k=n\}$ is the
divisor function of order $k$. The generalized Ramanujan conjecture predicts that one may take
$\theta_{m}=0$ in \eqref{theta}.

Define the analytic conductor
\[
 \mathcal C(\pi,t)= q_\pi
 \prod_{j=1}^m
 \bigl(3+|it+\nu_{\pi,j}|\bigr),
 \qquad \mathcal C(\pi)=\mathcal C(\pi,0).
\]
By applying
the Phragm\'en--Lindel\"of principle and the functional equation, we have
the convexity bound (see \cite[(5.20)]{IK04}):
\begin{equation*}
 |L(\sigma+it,\pi)|
 \ll\mathcal C(\pi,t)^{(1-\sigma)/2+\varepsilon}, \qquad \quad
 0\le\sigma\le1
\end{equation*}
for any $\varepsilon>0$.

\subsection{Rankin-Selberg $L$-functions}

Let $\pi=\otimes_v\pi_{v}\in \mathcal{A}(m)$
and $\pi'=\otimes_v\pi'_{v}\in \mathcal{A}(m')$.
The
Rankin-Selberg \(L\)-function associated
to $\pi$ and $\pi'$ is of the form
\[
 L(s,\pi\times\pi')
 =\prod_p L(s,\pi_{p}\times\pi'_p)
 =\sum_{n=1}^{\infty}
   \frac{\lambda_{\pi\times\pi'}(n)}{n^s},
\]
where the Euler product and Dirichlet series converge absolutely for
\(\mathrm{Re}(s)>1\), and for each finite prime $p$
\bea\label{rankin-p}
 L(s,\pi_{p}\times\pi'_{p})^{-1}
 =\prod_{j=1}^{m}\prod_{j'=1}^{m'}
  \left(1-\alpha_{\pi\times\pi',j,j'}(p)p^{-s}\right)=\sum_{k=0}^{mm'}
\frac{\mu_{\pi\times \pi'}(p^k)}{p^{ks}}
\eea
for suitable complex numbers $\alpha_{\pi\times\pi',j,j'}(p)$.
In particular, some of the local parameters $\alpha_{\pi\times\pi',j,j'}(p)$ are zero if the degree of
$L(s,\pi_{p}\times\pi'_{p})^{-1}$ is less than $mm'$.
At the infinity place there exist $mm'$ complex Langlands parameters
$\nu_{\pi\times\pi',j,j'}$ associated to $\pi_{\infty}$ and $\pi'_{\infty}$, from which one defines
\[
 L(s,\pi_{\infty}\times \pi'_{\infty}):=\prod_{j=1}^m\prod_{j'=1}^{m'}
 \Gamma_{\mathbb R}(s+\nu_{\pi\times\pi',j,j'}).
\]
The completed \(L\)-function
\[
 \Lambda(s,\pi\times\pi')
 =q_{\pi\times\pi'}^{s/2}
  L(s,\pi_{\infty}\times\pi'_{\infty})
  L(s,\pi\times\pi')
\]
has a meromorphic continuation to $\mathbb C$ and satisfies
the functional equation
\[
 \Lambda(s,\pi\times\pi')
 =\varepsilon(\pi\times\pi')
  \Lambda(1-s,\widetilde\pi\times\widetilde{\pi}'),
 \qquad
 |\varepsilon(\pi\times\pi')|=1,
\]
where \(q_{\pi\times\pi'}\) is the arithmetic conductor
of $\pi\times \pi'$.
The finite Rankin--Selberg $L$-function $L(s,\pi\times\pi')$ is entire unless
$m=m'$ and $\pi'\simeq\widetilde\pi\otimes|\det|^{it_0}$ for some $t_0\in\mathbb R$.
With the convention $L(s,\tau\otimes|\det|^{it_0})=L(s+it_0,\tau)$,
its pole on $\operatorname{Re}s=1$ is then the simple pole at $1-it_0$.
In particular, $L(s,\pi\times\widetilde\pi)$ has a simple pole at $1$ with positive residue.
We also use nonvanishing on $\operatorname{Re}s=1$ away from the pole
\cite{JacquetPiatetskiShapiroShalika1983,Brum}.

With $\theta_m$ as in \eqref{theta}, we have the pointwise bound
\[
 \bigl|\alpha_{\pi\times\pi',j,j'}(p)\bigr|
 \leq p^{\theta_{m}+\theta_{m'}},
 \quad\qquad
\mathrm{Re}(\nu_{\pi\times \pi',j,j'})\ge -\theta_{m}-\theta_{m'}
\]
and
\bea\label{pi-pi}
 \bigl|\lambda_{\pi\times\pi'}(n)\bigr|
 \leq d_{mm'}(n)n^{\theta_{m}+\theta_{m'}}.
\eea
Moreover, if \(p\nmid q_{\pi} q_{\pi'}\), we
have the equality of sets
\bea\label{unrami-p}
 \left\{\alpha_{\pi\times\pi',j,j'}(p):
 1\leq j\leq m,1\leq j'\leq m'\right\}
 =
 \left\{
  \alpha_{\pi,j}(p)\alpha_{\pi',j'}(p):
  1\leq j\leq m,\ 1\leq j'\leq m'
 \right\}.
\eea

Since the unramified parameters of \(\widetilde\pi\) are the complex
conjugates of those of \(\pi\), for $p\nmid q_\pi$, one has
\[
 \lambda_{\pi\times\widetilde\pi}(p)
   =|\lambda_{\pi}(p)|^2.
\]
More generally, the coefficients of
\(L(s,\pi\times\widetilde\pi)\) are nonnegative, and the following
inequalities hold for every \(n\geq1\)
(see \cite[(3.4), (3.5), and Lemma~3.1]{J25}):
\begin{equation}\label{ineq:Coefficient&Reciprocal}
 \begin{aligned}
 |\lambda_{\pi\times\pi'}(n)|
 &\leq
 \bigl(
   \lambda_{\pi\times\widetilde\pi}(n)
   \lambda_{\pi'\times\widetilde \pi'}(n)
 \bigr)^{1/2},\\
 |\mu_{\pi\times\pi'}(n)|
 &\leq
 \bigl(
   \lambda_{\pi\times\widetilde\pi}(n)
   \lambda_{\pi'\times\widetilde\pi'}(n)
 \bigr)^{1/2}.
 \end{aligned}
\end{equation}
Therefore,
\bea\label{mu-upper}
 |\lambda_{\pi}(n)|^2
 \leq\lambda_{\pi\times\widetilde\pi}(n),
 \qquad
 |\mu_{\pi}(n)|^2
 \leq\lambda_{\pi\times\widetilde\pi}(n),
 \qquad
 |\mu_{\pi\times\widetilde\pi}(n)|
 \leq\lambda_{\pi\times\widetilde\pi}(n).
\eea

Define the analytic conductor of $\pi\times \pi'$ by
\[
 \mathcal C(\pi\times \pi',t)= q_{\pi\times \pi'}
 \prod_{j=1}^m\prod_{j'=1}^{m'}
 \bigl(3+|it+\nu_{\pi\times\pi',j,j'}|\bigr),
 \qquad \mathcal C(\pi\times \pi')=\mathcal C(\pi\times \pi',0).
\]
Then $\mathcal C(\pi\times \pi',t)$ satisfies the inequalities (see \cite[(5.8), (5.11)]{IK04}):
\[
 \mathcal C(\pi\times \pi',t)
 \ll C(\pi\times \pi')(1+|t|)^{mm'},
 \qquad
 \mathcal C(\pi\times\pi')
 \ll_{m,m'}
 \mathcal C(\pi)^{m'}\mathcal C(\pi')^m.
\]
Write \(s=\sigma+it\), where \(0\leq\sigma\leq1\). If
\(\pi'\simeq\widetilde\pi\otimes|\det|^{it_0}\), put
\(b_{\pi,\pi'}(s)=(s+it_0-1)/(s+it_0-2)\). In all other cases put
\(b_{\pi,\pi'}(s)=1\). For every \(\varepsilon>0\), the convexity bound,
with the pole removed, takes the form
\begin{equation*}
 \bigl|b_{\pi,\pi'}(s)L(s,\pi\times\pi')\bigr|
 \ll_{\pi,\pi',\varepsilon}
 \mathcal C(\pi\times\pi',t)^{(1-\sigma)/2+\varepsilon}.
\end{equation*}
At \(s=1-it_0\) the product is understood by continuity.

In many applications, we use the following important estimate of Jiang \cite{J25}.

\begin{theorem}\label{Jiang}
Let $\pi\in\mathcal{A}(m)$.  For any $\varepsilon>0$ and any
$\sigma\geq1-1/(m^2+1)+\varepsilon$,
\[
 \prod_p\left(1+\sum_{k=2}^\infty
 \frac{\lambda_{\pi\times\widetilde\pi}(p^k)}{p^{k\sigma}}\right)
 \ll_{m,\varepsilon}\mathcal C(\pi)^\varepsilon.
\]
\end{theorem}

Let $\chi\bmod q$ be induced by the primitive character $\chi^*$.
In the notation of \cite[Section~7]{J25}, we have for $\operatorname{Re}(s)>1$
\begin{align*}
 L(s,\pi\otimes\chi)
 &:=\sum_{n\ge1}\frac{\lambda_\pi(n)\chi(n)}{n^s}
   =L(s,\pi\times\chi^*)P_{\pi,\chi;q}(s),\\
 P_{\pi,\chi;q}(s)&:=\prod_{p\mid q}L(s,\pi_p\times\chi_p^*)^{-1}.
\end{align*}
For primitive $\chi\bmod q$, the conductor and analytic conductor satisfy
\[
 q_{\pi\times\chi}\ll_\pi q^m,\qquad
 \mathcal C(\pi\times\chi,t)\ll_\pi\bigl(q(3+|t|)\bigr)^m.
\]

\subsection{Isobaric $L$-functions}

We  return to the fixed representations of Section~1, with $m_j\ge2$ and
$\pi_h\not\simeq\pi_j\otimes|\det|^{it}$ for $h\ne j$ and $t\in\mathbb R$.
With a slight abuse of notation, throughout the rest of the paper,
we reuse the symbol
$\pi$ to denote the isobaric representation
\[
\pi=\boxplus_{j=1}^r\pi_j=\pi_1\boxplus\cdots\boxplus\pi_r
\]
of $\mathrm{GL}_m(\mathbb A_{\mathbb Q})$, where $m=\sum_{j=1}^rm_j$,
which should not be confused with the one defined in Subsections 2.1--2.2.

The global $L$-function attached to $\pi$ is of degree $m$ and
factors as a product of the $L$-functions of its constituent representations:
\bea\label{factor}
 L(s,\pi)=\prod_{j=1}^rL(s,\pi_j).
\eea
For $\mathrm{Re}(s)>1$, define
\bna
 L(s,\pi)=\sum_{n=1}^{\infty}\frac{\lambda_{\pi}(n)}{n^s}, \qquad  L(s,\pi_j)=\sum_{n=1}^{\infty}\frac{\lambda_{\pi_j}(n)}{n^s}.
\ena
Then by \eqref{factor},
\bea\label{Fourier coeffi-factor}
\lambda_{\pi}(n)=\sum_{n=n_1n_2\cdots n_r}
\lambda_{\pi_1}(n_1)\lambda_{\pi_2}(n_2)\cdots\lambda_{\pi_r}(n_r):=
\lambda_{\pi_1}*\lambda_{\pi_2}*\cdots *\lambda_{\pi_r}(n).
\eea

The associated Rankin-Selberg $L$-function
is defined by
\bea\label{Rankin-Selberg def}
L(s,\pi \times \widetilde{\pi})=\prod_{h,j=1}^r L(s,\pi_h \times \widetilde\pi_j),
\eea
for $\mathrm{Re}(s)\gg 1$. The pole at $s=1$ has order $r$; the real-twist inequivalence assumption
excludes off-diagonal translated poles. By
Lemma a in Hoffstein and Ramakrishnan \cite{HR}, $\lambda_{\pi\times \widetilde{\pi}}(n)\geq 0$.
For $\mathrm{Re}(s)\gg 1$, define
\bea\label{Rankin-Selberg-1}
L(s,\pi \times \widetilde{\pi})^{-1}&=&\prod_{h,j=1}^r L(s,\pi_h \times \widetilde\pi_j)^{-1}:=\sum_{n=1}^{\infty}
\frac{\mu_{\pi\times \widetilde{\pi}}(n)}{n^s},\\
L(s,\pi_h \times \widetilde{\pi}_j)^{-1}&=&
\sum_{n=1}^{\infty}
\frac{\mu_{\pi_h\times \widetilde{\pi}_j}(n)}{n^s}\nonumber
\eea
Then by \eqref{Rankin-Selberg def},
\bea\label{Fourier coeffi-factor-mu}
\mu_{\pi\times \widetilde{\pi}}(n)&=&\sum_{n=n_{11}n_{12}\cdots n_{rr}}
\prod_{h,j=1}^r
\mu_{\pi_h\times \widetilde{\pi}_j}(n_{hj})\nonumber\\
&:=&
\mu_{\pi_1\times \widetilde{\pi}_1}*\mu_{\pi_1\times \widetilde{\pi}_2}
*\cdots *\mu_{\pi_r\times \widetilde{\pi}_r}(n).
\eea
Here, in particular, for primes \(p\) at which all the \(\pi_j\) are unramified, by \eqref{rankin-p} and \eqref{unrami-p},
\[
L(s,\pi_{h,p} \times \widetilde\pi_{j,p})^{-1}
=
\prod_{\ell=1}^{m_h}\prod_{\ell'=1}^{m_j}
\left(1-\alpha_{\pi_h,\ell}(p)\overline{\alpha_{\pi_j,\ell'}(p)}p^{-s}\right)
=
\sum_{k=0}^{m_h m_j} \mu_{\pi_h\times \widetilde{\pi}_j}(p^k) p^{-ks},
\]
where the coefficients $\mu_{\pi_h\times \widetilde{\pi}_j}(p^k)$ are given by elementary symmetric functions:
\[
\mu_{\pi_h\times \widetilde{\pi}_j}(p^k)
=
(-1)^k
\sum_{\substack{S\subseteq \{1,\dots,m_h\}\times\{1,\dots,m_j\}\\ |S|=k}}
\prod_{(\ell,\ell')\in S}\alpha_{\pi_h,\ell}(p)\overline{\alpha_{\pi_j,\ell'}(p)}.
\]

For any Dirichlet character $\chi$, complete multiplicativity gives
\[
 L(s,\pi\otimes\chi)=\prod_{j=1}^rL(s,\pi_j\otimes\chi).
\]
The comparison with $L(s,\pi\times\chi^*)$ in Section~2.2 thus also
holds for the isobaric sum. For primitive $\chi$ of conductor $q$,
twisting each constituent gives the isobaric representation $\Pi_\chi$,
with standard $L$-function and analytic conductor
\[
 L(s,\Pi_\chi)=L(s,\pi\times\chi)=\prod_{j=1}^rL(s,\pi_j\times\chi),
 \qquad
 \mathcal C(\Pi_\chi,t)=\prod_{j=1}^r\mathcal C(\pi_j\times\chi,t).
\]
In particular,
\begin{align}\label{eq:isobaric-L-degree&conductor}
 q_{\Pi_\chi}\ll_\pi q^m,\qquad
 \mathcal C(\Pi_\chi,t)\ll_\pi\bigl(q(3+|t|)\bigr)^m .
\end{align}
The conductor bound of Bushnell and Henniart \cite{BushnellHenniart}
also gives $q_{\Pi_\chi}\le q_\pi q^m/(q_\pi,q)$, where $q_\pi$ is
the arithmetic conductor of $\pi$.

\subsection{Mean square estimates of \(L\)-functions}

We quote the following general theorem on mean values of analytic
functions as in \cite[Section~7.8]{T86}.
\begin{lem}[Mean-square convexity lemma]\label{lem:mean-value-interpolation}
Let $f(s)$ be an analytic function of $s$, real for real $s$, regular for $\sigma \ge \alpha$ except possibly
for a pole at $s = s_0$, and $O(e^{\epsilon|t|})$ as $|t| \to \infty$ for every positive $\epsilon$
and $\sigma \ge \alpha$. Let $\alpha < \beta$, and suppose that for all $U\geq 1$
\[
\int_{0}^{U} |f(\alpha+it)|^{2} \dif t \le CU^{a},
\]
\[
\int_{0}^{U} |f(\beta+it)|^{2} \dif t \le C'U^{b},
\]
where $a \ge 0$, $b \ge 0$, and $C$, $C'$ depend on $f(s)$. Then for $\alpha < \sigma < \beta$, $U \ge 2$,
\[
\int_{U}^{2U} |f(\sigma+it)|^{2} \dif t \ll (CU^{a})^{(\beta-\sigma)/(\beta-\alpha)}(C'U^{b})^{(\sigma-\alpha)/(\beta-\alpha)},
\]
where the implied constant depends on $a,b,\alpha,\beta$ only.
\end{lem}

The following mean-square estimate applies in particular to the isobaric
representation $\pi$.
\begin{lem}\label{meanLbound}
Let $\mathcal F$ be associated with a fixed entire automorphic $L$-function of degree $D\geq2$.
For $1/2\leq\sigma\leq1+\varepsilon$ and $U\geq 2$,
\[
 \int_U^{2U}|L(\sigma+it,\mathcal F)|^2\dif t
 \ll_{\mathcal F,\varepsilon}
 U^{D-1+\sigma(2-D)+D\varepsilon}.
\]
\end{lem}
\begin{proof}
On the critical line $\Re(s)=1/2$, the approximate functional equation
and large sieve estimates from \cite[Theorems~5.3 and 7.17]{IK04} give
\[
 \int_U^{2U}|L(1/2+it,\mathcal F)|^2\dif t
 \ll_{\mathcal F,\varepsilon}U^{D/2+\varepsilon}.
\]
On the other hand, absolute convergence gives a mean square $\ll_{\mathcal F,\varepsilon}U\le U^{1+D\varepsilon}$ on $\Re(s)=1+\varepsilon$.
Let \(F(s)=L(s,\mathcal F)\) and
\(F^\#(s)=\overline{F(\bar s)}=L(s,\widetilde{\mathcal F})\). The same
endpoint estimates hold for \(F^\#\), since the contragredient has the
same degree. The entire functions
\[
 F_1=(F+F^\#)/2,\qquad F_2=(F-F^\#)/(2i)
\]
are real on the real axis and satisfy both endpoint bounds. Summing the
dyadic bounds gives the corresponding estimates on \([0,U]\).
Applying Lemma~\ref{lem:mean-value-interpolation} to \(F_1\) and \(F_2\)
and using \(|F|^2\leq2(|F_1|^2+|F_2|^2)\) gives the claimed exponent.
\end{proof}

\section{Proof of Theorems~\ref{mainresult1} and \ref{mainresult2}}

Let $V_\pi(Q,X)$ be defined as in \eqref{DeltaPiQX}. We write
\bea\label{V-total}
 V_\pi(Q,X)=\Delta_\pi(Q,X)-C_\pi(Q,X),
 \eea
 where
 \bea\label{C}
 C_\pi(Q,X)=\sum_q\frac{\Phi(q/Q)}{\phi(q)}
 \left|\sum_{(n,q)=1}\lambda_\pi(n)\Psi\!\left(\frac nX\right)\right|^2,
\eea
and
\[
 \Delta_\pi(Q,X)=
 \sum_{q\geq1}\Phi\!\left(\frac qQ\right)
 \sum_{\substack{a\!\!\!\!\! \!\pmod q\\(a,q)=1}}\;
 \Bigl| \sum_{n\equiv a\!\!\!\!\! \!\pmod q}
 \lambda_\pi(n)\Psi\!\left(\frac nX\right)\Bigr|^2.
\]

We first show that $C_\pi(Q,X)$ is small. In fact,
Mellin inversion gives
\[
\sum_{(n,q)=1}\lambda_\pi(n)\Psi\!\left(\frac nX\right)=\frac{1}{2\pi i}\int_{(2)}\widetilde\Psi(s)X^sL(s,\pi)
 \prod_{p\mid q}L(s,\pi_p)^{-1}\,\mathrm{d}s.
\]
Since $m_j\ge2$, $L(s,\pi)$ is entire. Each inverse local factor is a
polynomial. Moving the contour to $\operatorname{Re}(s)=1/2+\varepsilon$ and by rapid
Mellin decay, we have, for every $\varepsilon>0$,
$$\sum_{(n,q)=1}\lambda_\pi(n)\Psi\!\left(\frac nX\right)\ll_{\pi,\Psi}X^{1/2+\varepsilon}q^\varepsilon.$$
It follows that $C_\pi(Q,X)$ in \eqref{C} is bounded by
$$
C_\pi(Q,X)\ll X(QX)^{\varepsilon}.
$$
By \eqref{V-total} and the preceding bound for \(C_\pi(Q,X)\), it suffices
to estimate \(\Delta_\pi(Q,X)\) in the range of each theorem. Fix
\(\varepsilon>0\). For Theorem~\ref{mainresult1} the unconditional range is
\(Q^\varepsilon\leq X\leq Q^{1+2/(m+2)-\varepsilon}\).
For Theorem~\ref{mainresult2}, under the GRC and GRH assumptions stated
there, the range is \(Q^\varepsilon\leq X\leq Q^{2-\varepsilon}\).
In each case we shall prove, uniformly in \(X\) as \(Q\to\infty\),
\bna
 \Delta_\pi(Q,X)=QX P_{r-1}(\log X)+o(QX).
\ena

Expanding the square, we obtain
\bna
 \Delta_\pi(Q,X)= \mathscr{A}_\pi(Q,X)+ \mathscr{B}_\pi(Q,X),
\ena
where $\mathscr{A}_\pi(Q,X)$ and $\mathscr{B}_\pi(Q,X)$ denote the contributions arising from the diagonal
and off-diagonal terms in the expansion, respectively, defined by
\begin{align}
\mathscr{A}_\pi(Q,X)
 &=\sum_q\Phi\left(\frac{q}{Q}\right)
 \sum_{(n,q)=1}|\lambda_\pi(n)|^2
 \Psi\left(\frac{n}{X}\right)^2,\label{D}\\
\mathscr{B}_\pi(Q,X)
 &=\sum_q\Phi\left(\frac{q}{Q}\right)
   \sum_{\substack{n_1\ne n_2\\(n_1n_2,q)=1\\n_1\equiv n_2(\bmod q)}}
   \lambda_\pi(n_1)\overline{\lambda_\pi(n_2)} \Psi\left(\frac{n_1}{X}\right)
 \Psi\left(\frac{n_2}{X}\right).\label{OD}
\end{align}

First, we transform the off-diagonal contribution $\mathscr{B}_\pi(Q,X)$ in \eqref{OD}.
Write
\begin{align*}
 \mathscr B_\pi(Q,X)
 ={}&
 \sum_{\substack{n_1=dn_1', n_2=dn_2'\\(n_1',n_2')=1\\n_1'\ne n_2'}}
 \lambda_\pi(n_1)\overline{\lambda_{\pi}(n_2)}\Psi\left(\frac{n_1}{X}\right)
 \Psi\left(\frac{n_2}{X}\right)
 \sum_{\substack{q\mid (n_1'-n_2')\\(q, dn_1'n_2')=1}}
 \Phi\left(\frac qQ\right).
\end{align*}
Note that the conditions \((n_1',n_2')=1\) and \(q|(n_1'-n_2')\) imply
\((n_1'n_2',q)=1\). Then the innermost sum reduces to
\begin{align*}
 \sum_{\substack{q\mid (n_1'-n_2')\\(q,d)=1}}
 \Phi\left(\frac qQ\right)
 =
 \sum_{q\mid (n_1'-n_2')}\Phi\left(\frac qQ\right)
 \sum_{c\mid(q,d)}\mu(c)=
 \sum_{c\mid d}\mu(c)
 \sum_{c\ell\mid (n_1'-n_2')}
 \Phi\left(\frac{c\ell}{Q}\right).
\end{align*}
Let $|n_1'-n_2'|=c\ell g$. The condition \(c\ell| (n_1'-n_2')\) is equivalent to \(cg|(n_1'-n_2')\).
Hence
\begin{align}\label{eq:r-reparametrization}
 \sum_{\substack{q\mid (n_1'-n_2')\\(q,d)=1}}
 \Phi\left(\frac qQ\right)
 =
 \sum_{c\mid d}\mu(c)
 \sum_{\substack{g\geq1\\cg\mid (n_1'-n_2')}}
 \Phi\left(\frac{|n_1'-n_2'|}{gQ}\right).
\end{align}
Moreover, the conditions $(n_1',n_2')=1$ and \(n_1'\equiv n_2'\pmod{cg}\)
implies \((n_1'n_2',cg)=1\).
Detecting the condition \(cg|(n_1'-n_2')\) by
Dirichlet characters modulo \(cg\), one has
\[
 \mathbf 1_{cg\mid (n_1'-n_2')}
 =
 \frac{1}{\phi(cg)}
 \sum_{\chi(\bmod cg)}
 \chi(n_1')\overline{\chi(n_2')}.
\]
Substituting this into
\eqref{eq:r-reparametrization}, we obtain
\begin{align*}
 \sum_{\substack{q\mid (n_1'-n_2')\\(q,d)=1}}
 \Phi\left(\frac qQ\right)
 =
 \sum_{g\geq1}\sum_{c\mid d}
 \frac{\mu(c)}{\phi(cg)}
 \sum_{\chi(\bmod cg)}
 \chi(n_1')\overline{\chi(n_2')}
 \Phi\left(\frac{|n_1'-n_2'|}{gQ}\right).
\end{align*}
Therefore,
\bna
 \mathscr B_\pi(Q,X)
 =\sum_{\substack{n_1=dn_1'\\ n_2=dn_2'\\(n_1',n_2')=1}}
 \lambda_\pi(n_1)\overline{\lambda_{\pi}(n_2)}
 \sum_{g\geq1}\sum_{c\mid d}
 \frac{\mu(c)}{\phi(cg)}
 \sum_{\chi(\bmod cg)}\chi(n_1')\overline{\chi(n_2')}\,
 \mathcal W\left(\frac{n_1}{X},\frac{n_2}{X};\frac{X}{d g Q}\right),
\ena
where
\bna
 \mathcal W(x_1,x_2;y)=\Psi(x_1)\Psi(x_2)\Phi\bigl(|x_1-x_2|y\bigr).
\ena
Here the condition \(n_1'\ne n_2'\) is omitted, since
\(\Phi\in C_c^\infty((0,\infty))\) vanishes in a neighborhood of zero.

Let \(\chi_0\) denote the principal character modulo \(cg\). By
separating \(\chi_0\) from the nonprincipal characters, we have
\bna
 \mathscr B_\pi(Q,X)
 =\mathscr B_{1,\pi}(Q,X)+\mathscr B_{2,\pi}(Q,X),
\ena
where
\begin{align*}
\mathscr B_{1,\pi}(Q,X)
 ={}&\sum_{\substack{n_1=dn_1'\\ n_2=dn_2'\\(n_1',n_2')=1}}
 \lambda_\pi(n_1)\overline{\lambda_{\pi}(n_2)}
 \sum_{g\geq1}
 \sum_{\substack{c\mid d\\(n_1'n_2',cg)=1}}
 \frac{\mu(c)}{\phi(cg)}\,
 \mathcal W\left(\frac{n_1}{X},\frac{n_2}{X};\frac{X}{dgQ}\right),
 \\
 \mathscr B_{2,\pi}(Q,X)
 ={}&\sum_{\substack{n_1=dn_1'\\ n_2=dn_2'\\(n_1',n_2')=1}}
 \lambda_\pi(n_1)\overline{\lambda_{\pi}(n_2)}
 \sum_{g\geq1}\sum_{c\mid d}
 \frac{\mu(c)}{\phi(cg)}
 \sum_{\substack{\chi(\bmod cg)\\ \chi\ne\chi_0}}
 \chi(n_1')\overline{\chi(n_2')}\,
 \mathcal W\left(\frac{n_1}{X},\frac{n_2}{X};\frac{X}{dgQ}\right).
\end{align*}

Theorems~\ref{mainresult1} and \ref{mainresult2} follow from
Propositions~\ref{finald}, \ref{finalg1}, and \ref{finalg2}.

\begin{proposition}\label{finald}
Let \(\pi,\Phi,\Psi\) satisfy the hypotheses of Section~1, and fix
\(0<\delta<B\). There is a polynomial \(P_{r-1}\) of degree \(r-1\)
such that, as \(Q\to\infty\), uniformly for \(Q^\delta\leq X\leq Q^B\),
\[
\mathscr{A}_\pi(Q,X)=QX P_{r-1}(\log X)+o(QX),
\]
where the leading coefficient is the constant given by
\begin{align*}
 \frac{B_{\pi,0}(1)H_\pi(1)}{(r-1)!}
 \prod_{h=1}^r
 \operatorname*{Res}_{s=1}
 L(s,\pi_h\times\widetilde\pi_h)
 \prod_{\substack{h,j=1 \\ h\neq j}}^r
 L(1,\pi_h\times\widetilde\pi_j)
\end{align*}
with
\begin{align*}
B_{\pi,0}(1)=\prod_p\bigg\{1+
\frac{1}{p}
 \bigg(\Big(1+\sum_{\ell\geq1}\frac{|\lambda_{\pi}(p^\ell)|^2}{p^\ell}\Big)^{-1}-1\bigg)
 \bigg\}.
\end{align*}
\end{proposition}

\begin{proposition}\label{finalg1}
Let $\pi$ satisfy the hypotheses of Section~1.
If $m\leq2$, then unconditionally
$ \mathscr B_{1,\pi}(Q,X)=o(QX)$ for $X\leq Q^{2-\varepsilon}$.
If $m>2$, then unconditionally
\[
  \mathscr B_{1,\pi}(Q,X)=o(QX)
 \quad\text{for}\quad
 X\leq Q^{1+2/(m+2)-\varepsilon}.
\]
Under GRH and GRC the estimate holds for $X\leq Q^{2-\varepsilon}$.
\end{proposition}

\begin{proposition}\label{finalg2}
Let $\pi$ satisfy the hypotheses of Section~1.
Unconditionally,
\[
 \mathscr B_{2,\pi}(Q,X)
 \ll_{\pi,\Phi,\Psi,\varepsilon}
 X^{1+\varepsilon}
 \left\{\frac XQ+\left(\frac XQ\right)^{m/2}\right\}.
\]
Consequently $ \mathscr B_{2,\pi}(Q,X)=o(QX)$ for
$X\leq Q^{1+\min\{1,2/m\}-\varepsilon}$.
Under GRH the estimate holds for $X\leq Q^{2-\varepsilon}$.
\end{proposition}

We defer their proofs to
Sections~\ref{sec:proof-finald}, \ref{sec:proof-finalg1}, and
\ref{sec:proof-finalg2}, respectively.

\section{The diagonal contribution $\mathscr{A}_\pi(Q,X)$}
\label{sec:proof-finald}

In this section, we prove Proposition~\ref{finald}. Recall the definition of $\mathscr{A}_\pi(Q,X)$ in
\eqref{D}, which we relabel as
\begin{align}
\mathscr{A}_\pi(Q,X)
 =\sum_q\Phi\left(\frac{q}{Q}\right)
 \sum_{(n,q)=1}|\lambda_\pi(n)|^2
 \Psi\left(\frac{n}{X}\right)^2,\label{D-}
\end{align}
where $\lambda_\pi(n)$ are the Dirichlet coefficients of the
isobaric $L$-function $L(s,\pi)$,
$\pi=\pi_1\boxplus\cdots\boxplus\pi_r$ with $\pi_j$
unitary cuspidal automorphic
representations of
$\mathrm{GL}_{m_j}(\mathbb{A}_{\mathbb{Q}})$, $j=1,\ldots,r$,
and $m=\sum_{j=1}^rm_j$.
Here $\Phi,\Psi\in C_c^\infty((0,\infty);\mathbb R)$, $\Phi\ge0$, and
    $\displaystyle\int_0^\infty \Phi(x)\,\dif x=
    \int_0^\infty \Psi(x)^2\,\dif x=1$.

Denote by $\widetilde{\Psi^2}(s)$ the Mellin transform of $\Psi(x)^2$.
By the Mellin inversion formula, one has
\[
\sum_{(n,q)=1}|\lambda_\pi(n)|^2\Psi\!\left(\frac nX\right)^2
=
\frac{1}{2\pi i}\int_{(\sigma)}
\widetilde{\Psi^2}(s)X^s
\sum_{(n,q)=1}\frac{|\lambda_\pi(n)|^2}{n^s}\,\dif s,
\]
where $\sigma>1$. Note that for $\mathrm{Re}(s)>1$, one has the Euler product
\bea\label{D_p}
\sum_{n=1}^\infty \frac{|\lambda_\pi(n)|^2}{n^s}
=\prod_p\bigg(\sum_{k\ge 0}\frac{|\lambda_\pi(p^k)|^2}{p^{ks}}\bigg)
=:\prod_p D_p(s),
\eea
say. Hence
\begin{equation}
\sum_{(n,q)=1}|\lambda_\pi(n)|^2\Psi\!\left(\frac nX\right)^2
=
\frac{1}{2\pi i}\int_{(\sigma)}
\widetilde{\Psi^2}(s)X^s
\bigg(\sum_{n=1}^\infty \frac{|\lambda_\pi(n)|^2}{n^s}\bigg)
\prod_{p\mid q}D_p(s)^{-1}\,\dif s.
\label{eq:5.2.2}
\end{equation}

Moreover, we obtain the following result, whose proof
is deferred to the end of this section.
 \begin{lem}\label{corollary:Jiang}
Let $m^*=\max_{1\leq j\leq r}m_j$.
For \(\Re s>1\) one has
\begin{align}
 \sum_{n=1}^\infty\frac{|\lambda_\pi(n)|^2}{n^s}
 &=H_\pi(s)\prod_{h,j=1}^r
   L(s,\pi_h\times\widetilde\pi_j),\label{eq:H-factor}\\
 \prod_p\left(1-\frac{|\lambda_{\pi}(p)|^2}{p^s}\right)
 &=H_{\pi}^{\sharp}(s)\prod_{h,j=1}^r
   L(s,\pi_h\times\widetilde\pi_j)^{-1},\label{eq:H1-factor}
\end{align}
The Euler products defining \(H_\pi\) and \(H_\pi^\sharp\) converge
absolutely and locally uniformly for \(\Re s>1-1/(m^{*2}+1)\).
They are holomorphic there and satisfy
\[
 |H_\pi(s)|+|H_\pi^\sharp(s)|\ll_{\pi,\varepsilon}\mathcal C(\pi)^\varepsilon
 \qquad\bigl(\Re s\geq 1-1/(m^{*2}+1)+\varepsilon\bigr).
\]
Under GRC the two correction products converge absolutely and locally
uniformly throughout \(\Re s>1/2\) and are holomorphic there. In each
case the right-hand sides of \eqref{eq:H-factor} and \eqref{eq:H1-factor}
give meromorphic continuations to the stated half-plane.
\end{lem}

From \eqref{D-}, \eqref{eq:5.2.2} and \eqref{eq:H-factor}, we obtain
\begin{equation}
\mathscr A_\pi(Q,X)
=
\frac{1}{2\pi i}\int_{(\sigma)}
\widetilde{\Psi^2}(s)X^s
H_\pi(s) \prod_{h,j=1}^r
 L(s,\pi_h\times\widetilde\pi_j)
\sum_{q\ge1}\Phi\!\left(\frac qQ\right)\prod_{p\mid q}D_p(s)^{-1}\,ds,
\label{eq:4.6}
\end{equation}
where we have interchanged summation and integration.

Next we shift the contour of integration in \eqref{eq:4.6} to
$\operatorname{Re}(s)=\sigma_0:=1-\frac{1}{m^{*2}+1}+\varepsilon$.
The only pole crossed in this contour shift is \(s=1\), which is of order \(r\). This is because each
\(L(s,\pi_h\times\widetilde\pi_h)\) has a simple pole at \(s=1\), $h=1,\dots,r$;
for \(h\ne j\), the exclusion of all unitary twists in Section~1 makes
\(L(s,\pi_h\times\widetilde\pi_j)\) entire, and nonvanishing on
\(\Re s=1\) gives a nonzero value at \(1\).
For each \(p\mid q\) the apparent singularities of \(D_p(s)^{-1}\)
cancel since \(H_{\pi,p}(s)D_p(s)^{-1}=L_p(s,\pi\times\widetilde\pi)^{-1}\).
At \(s=1\) every \(D_p(1)\) is positive, so absolute convergence of
the correction product gives \(H_\pi(1)\ne0\). Moreover, the horizontal integrals tend to zero owing to the
rapid decay of \(\widetilde{\Psi^2}(s)\) and the polynomial growth of the
Rankin--Selberg factors in vertical strips. Therefore, we have
\bea
\mathscr A_\pi(Q,X)
&=&\operatorname{Res}_{s=1}\Big\{
\widetilde{\Psi^2}(s)X^s
H_\pi(s) \prod_{h,j=1}^r
 L(s,\pi_h\times\widetilde\pi_j)
\sum_{q\ge1}\Phi\!\left(\frac qQ\right)\prod_{p\mid q}D_p(s)^{-1}
\Big\}\nonumber\\
&&+
\frac{1}{2\pi i}\int_{(\sigma_0)}
\widetilde{\Psi^2}(s)X^s
H_\pi(s) \prod_{h,j=1}^r
 L(s,\pi_h\times\widetilde\pi_j)
\sum_{q\ge1}\Phi\!\left(\frac qQ\right)\prod_{p\mid q}D_p(s)^{-1}\,ds,
\label{eq:4.7}
\eea

To compute its residue, set
\begin{align*}
 G_h(s)
 &=
 (s-1)L(s,\pi_h\times\widetilde\pi_h),\\
 A_\pi(s)
 &=
 \widetilde{\Psi^2}(s)H_\pi(s)
 \prod_{h=1}^rG_h(s)
 \prod_{\substack{h,j=1\\h\ne j}}^r
 L(s,\pi_h\times\widetilde\pi_j),\\
 V_q(s)
 &=
 \prod_{p\mid q}D_p(s)^{-1},\\
 U_q(s)
 &=
 A_\pi(s)V_q(s).
\end{align*}
The function \(U_q(s)\) is holomorphic in a neighborhood of \(s=1\).
Moreover, \(U_q(1)\ne0\), by the preceding nonvanishing statements,
Lemma~\ref{corollary:Jiang}, and the fact that \(D_p(1)>0\).
Write
\[
 U_q(s)
 =
 \sum_{\ell\geq0}u_{q,\ell}(s-1)^\ell,
 \qquad
 u_{q,\ell}
 =
 \frac{U_q^{(\ell)}(1)}{\ell!}.
\]
Since
\[
 X^s
 =
 X\sum_{k\geq0}
 \frac{(s-1)^k(\log X)^k}{k!},
\]
the residue at \(s=1\) is
\[
 XP_{q,r-1}(\log X),
\]
where
\begin{align}\label{eq:diagonal-polynomial}
 P_{q,r-1}(Y)
 =
 \sum_{\ell=0}^{r-1}
 \frac{u_{q,\ell}}{(r-1-\ell)!}
 Y^{r-1-\ell}.
\end{align}
In particular, \(P_{q,r-1}\) has degree \(r-1\), and its leading
coefficient is
\begin{align}\label{eq:q-leading-coefficient}
 \frac{H_\pi(1)}{(r-1)!}
 \prod_{h=1}^r
 \operatorname*{Res}_{s=1}
 L(s,\pi_h\times\widetilde\pi_h)
 \prod_{\substack{h,j=1\\h\ne j}}^r
 L(1,\pi_h\times\widetilde\pi_j)
 \prod_{p\mid q}D_p(1)^{-1},
\end{align}
where we used the normalization
\(\widetilde{\Psi^2}(1)=1\).

We next average all the coefficients in
\eqref{eq:diagonal-polynomial} over \(q\).  We begin with the leading
coefficient.  Since its dependence on \(q\) is contained in \(V_q(1)\),
we have
\begin{align*}
 \sum_{q=1}^{\infty}\frac{V_q(1)}{q^z}
 &=
 \prod_p
 \left(
  1+D_p(1)^{-1}\sum_{\ell\geq1}p^{-\ell z}
 \right)\notag\\
 &=
 \zeta(z)B_{\pi,0}(z),
\end{align*}
where
\[
 B_{\pi,0}(z)
 =
 \prod_p
 \left\{
  1+\bigl(D_p(1)^{-1}-1\bigr)p^{-z}
 \right\}.
\]
The local coefficient bounds give
\[
 D_p(1)^{-1}-1
 \ll_\pi p^{-1+2\theta},
\]
so the Euler product \(B_{\pi,0}(z)\) converges absolutely for
\(\mathrm{Re}(z)>2\theta\).  Mellin inversion against \(\Phi(q/Q)\), followed by
a shift of the \(z\)-contour to
\(\mathrm{Re}(z)=2\theta+\varepsilon\), therefore gives
\begin{align}\label{eq:average-Vq-zero}
 \sum_q\Phi(q/Q)V_q(1)
 =
 QB_{\pi,0}(1)+
 O_{\pi,\Phi,\varepsilon}
 \left(Q^{2\theta+\varepsilon}\right),
\end{align}
where we used \(\widetilde\Phi(1)=1\).

\bigskip

We now return to the proof of Lemma \ref{corollary:Jiang}.

\begin{proof}[Proof of Lemma \ref{corollary:Jiang}]
We first prove \eqref{eq:H-factor}. Write $q_{\pi}$ for the arithmetic conductor of $\pi$. Then
$q_{\pi}=\prod_{j=1}^r q_{\pi_j}$, where $q_{\pi_j}$ is the
arithmetic conductor of $\pi_j$.
For $p\nmid q_{\pi}$, by \eqref{Fourier coeffi-factor},
\[
 \lambda_{\pi}(p)=\sum_{h=1}^r\lambda_{\pi_h}(p),
\]
which implies
\bea\label{p term}
|\lambda_\pi(p)|^2
=
\Big|\sum_{h=1}^r \lambda_{\pi_h}(p)\Big|^2
=
\sum_{h,j=1}^r \lambda_{\pi_h}(p)\overline{\lambda_{\pi_j}(p)}
=
\sum_{h,j=1}^r \lambda_{\pi_h\times\widetilde\pi_j}(p).
\eea
This shows that the coefficients of \(p^{-s}\) agree in the Euler
products for
\[
\sum_{n=1}^\infty \frac{|\lambda_\pi(n)|^2}{n^s}
 \quad\text{and}\quad
 \prod_{h,j=1}^rL(s,\pi_h\times\widetilde\pi_j).
\]
Therefore, we define
\bea\label{H_p}
H_\pi(s):=\sum_{n=1}^\infty\frac{|\lambda_\pi(n)|^2}{n^s}\prod_{h,j=1}^r
   L(s,\pi_h\times\widetilde\pi_j)^{-1}:=\prod_p\Big(1+
       \sum_{k\geq 1}\frac{h_{\pi}(k,p)}{p^{ks}}\Big).
\eea

In the following, we estimate $h_{\pi}(k,p)$. By \eqref{D_p}, the local factor of $\sum_{n=1}^\infty |\lambda_\pi(n)|^2 n^{-s}$ at $p$ is
\[
1+\sum_{k\geq 1}|\lambda_{\pi}(p^k)|^2p^{-ks}.
\]
By \eqref{Rankin-Selberg-1}, the corresponding local factor of $\prod_{h,j=1}^r L(s,\pi_h\times\widetilde\pi_j)^{-1}$ at $p$ is
\begin{align*}
\sum_{k=0}^{m^2}\mu_{\pi\times \widetilde{\pi}}(p^k)p^{-ks}.
\end{align*}
Thus
\bea\label{eq:Hp-local-expansion}
 h_{\pi}(k,p)=\sum_{v=0}^k|\lambda_{\pi}(p^v)|^2\mu_{\pi\times \widetilde{\pi}}(p^{k-v}).
\eea
For $k=1$ and $p\nmid q_\pi$, by \eqref{p term},
\begin{align}\label{eq:linear-cancellation}
  h_{\pi}(1,p)= |\lambda_{\pi}(p)|^2+\mu_{\pi\times \widetilde{\pi}}(p)
   =\sum_{h,j=1}^r
     \lambda_{\pi_h\times\widetilde\pi_j}(p)+\mu_{\pi\times \widetilde{\pi}}(p)
   =0.
\end{align}

Set
\bea\label{max theta}
\theta=\frac{1}{2}-\frac{1}{m^{*2}+1}\ge\max_{1\le j\le r}\theta_{m_j},
\qquad m^*=\max_{1\le j\le r}m_j,
\eea
where $\theta_{j}$ is defined in \eqref{theta-j}.
By \eqref{Fourier coeffi-factor}, for $v\geq 1$,
\begin{align*}
 \lambda_{\pi}(p^v)=\sum_{k_1+\cdots+k_r=v}
       \prod_{j=1}^r\lambda_{\pi_j}(p^{k_j}),
\end{align*}
where by \eqref{lambda-n},
\[
 |\lambda_{\pi_j}(p^{k_j})|
 \le {k_j+m_j-1\choose m_j-1}p^{\theta_{m_j}k_j}\ll_{m_j} (k_j+1)^{m_j-1} p^{\theta k_j}.
\]
Here we have used the estimate ${u+d-1\choose d-1}\ll_d (u+1)^{d-1}$.
Order the tuples $(k_1, \dots, k_r)$ lexicographically and isolate the first nonzero exponent,
say $k_{j_0} = \ell$ (with $1 \leq \ell \leq v$). The remaining exponents sum to $v-\ell$.
Thus the product of the remaining factors
is bounded by $p^{\theta (v-\ell)}$ up to a polynomial
factor bounded by $O_m(v^{m-r})$. Moreover, the number of weak compositions \(k_1+\cdots+k_r=v\) is \(\binom{v+r-1}{r-1}\ll_r v^{r-1}\).
Consequently,
\bea\label{lambda-p-k-1}
|\lambda_\pi(p^v)| \ll_{m,r} v^{m-1} \sum_{j=1}^r \sum_{\ell=1}^v |\lambda_{\pi_j}(p^\ell)| p^{\theta(v-\ell)}.
\eea

Squaring both sides of \eqref{lambda-p-k-1} and applying the Cauchy--Schwarz inequality yields
\[
|\lambda_\pi(p^v)|^2 \ll_{m,r} v^{2m-1}  \sum_{j=1}^r
\sum_{\ell=1}^v |\lambda_{\pi_j}(p^\ell)|^2 p^{2\theta (v-\ell)}.
\]
Invoking the inequality $|\lambda_{\pi_j}(p^\ell)|^2 \leq \lambda_{\pi_j \times \widetilde{\pi}_j}(p^\ell)$ from \eqref{mu-upper}, we deduce
\bea\label{lambda-p-k-final}
|\lambda_\pi(p^v)|^2\ll_{m,r} v^{2m-1} \sum_{j=1}^r \sum_{\ell=1}^v \lambda_{\pi_j \times  \widetilde{\pi}_j}(p^\ell) p^{2\theta(v-\ell)}.
\eea

Similarly, we turn to estimating $\mu_{\pi\times \widetilde{\pi}}(p^k)$.
By \eqref{Fourier coeffi-factor-mu},
\bna
\mu_{\pi\times \widetilde{\pi}}(p^k)=\sum_{k_{hj}\geq 0\atop \sum_{h,j=1}^rk_{hj}=k}
\prod_{h,j=1}^r
\mu_{\pi_h\times \widetilde{\pi}_j}(p^{k_{hj}}).
\ena
Applying the second inequality in \eqref{ineq:Coefficient&Reciprocal} following with
the trivial bound $ab \ll (|a|+|b|)^2$, we obtain
\[
\left|\mu_{\pi_h\times \widetilde{\pi}_j}(p^{k_{hj}})\right|
\leq
\left(
\lambda_{\pi_h\times\widetilde{\pi}_h}(p^{k_{hj}})
\lambda_{\pi_j\times\widetilde{\pi}_j}(p^{k_{hj}})
\right)^{1/2}
\ll
\lambda_{\pi_h\times\widetilde{\pi}_h}(p^{k_{hj}})
+
\lambda_{\pi_j\times\widetilde{\pi}_j}(p^{k_{hj}}).
\]
Thus by \eqref{pi-pi},
\bna
\left|\mu_{\pi_h\times \widetilde{\pi}_j}(p^{k_{hj}})\right|
 \leq d_{m^{*2}}(p^{k_{hj}})p^{2\theta k_{hj}}\ll_{m}(k_{hj}+1)^{m^2-1}p^{2\theta k_{hj}}.
\ena
Order the pairs $(h,j)$ lexicographically
and isolate the first nonzero exponent,
say $k_{h_0 j_0} = \ell$ (with $1 \leq \ell \leq k$).
The remaining $r^2-1$ exponents sum to $k-\ell$,
and the product of the remaining factors
is bounded by $p^{2\theta (k-\ell)}$ up to a polynomial
factor in $k$ bounded by $O_m(k^{r^2m^2-r^2})$. Moreover, the number of weak compositions \(\sum_{h,j=1}^rk_{hj}=k\) is
\(\binom{k+r^2-1}{r^2-1}\ll_r k^{r^2-1}\).
Therefore
\bea\label{mu-pi-pi}
|\mu_{\pi \times \tilde{\pi}}(p^k)|
&\ll_{m}& k^{r^2m^2-1} \sum_{h,j=1}^r \sum_{\ell=1}^k \big|\mu_{\pi_h\times \widetilde{\pi}_j}(p^{\ell})\big| p^{2\theta(k-\ell)}\nonumber\\
&\ll_{m}& k^{r^2m^2-1} \sum_{h,j=1}^r \sum_{\ell=1}^k \Big(\lambda_{\pi_h\times\widetilde{\pi}_h}(p^{\ell})
+
\lambda_{\pi_j\times\widetilde{\pi}_j}(p^{\ell})\Big) p^{2\theta (k-\ell)}\nonumber\\
&\ll_{m,r}& k^{r^2m^2-1}
 \sum_{j=1}^r\sum_{\ell=1}^k
 \lambda_{\pi_j\times\widetilde\pi_j}(p^\ell)
 p^{2\theta (k-\ell)}.
\eea

Substituting \eqref{lambda-p-k-final} and \eqref{mu-pi-pi} into
\eqref{eq:Hp-local-expansion}, we obtain for $k\ge 2$,
\begin{align}\label{eq:h-local-majorant}
|h_{\pi}(k,p)| \ll_{m,r} k^{r^2 m^2 + 2m - 2}
&\biggl\{
\sum_{v=1}^{k-1}
\bigg(
\sum_{j=1}^r \sum_{\ell=1}^v \lambda_{\pi_j \times \widetilde{\pi}_j}(p^\ell) p^{2\theta(v-\ell)}
\bigg)\nonumber\\
&\quad\times
\bigg(
\sum_{j=1}^r \sum_{\ell=1}^{k-v} \lambda_{\pi_j \times \widetilde{\pi}_j}(p^\ell) p^{2\theta(k-v-\ell)}
\bigg)\nonumber\\
&\qquad+
\sum_{j=1}^r \sum_{\ell=1}^k
\lambda_{\pi_j \times \widetilde{\pi}_j}(p^\ell) p^{2\theta(k-\ell)}
\biggr\}\nonumber\\
=:k^{r^2 m^2 + 2m - 2}&
\biggl\{\sum_{v=1}^{k-1}B(v)B(k-v)+B(k)\biggr\},
\end{align}
where by \eqref{pi-pi},
\bna
B(v)&=&\sum_{j=1}^r \sum_{\ell=1}^v
\lambda_{\pi_j \times \widetilde{\pi}_j}(p^\ell) p^{2\theta(v-\ell)}\\
&\ll& \sum_{j=1}^r \sum_{\ell=1}^v
d_{m_j^2}(p^{\ell}) p^{2\theta_{m_j}\ell }p^{2\theta(v-\ell)}\\
&\ll_r& v^{m^2}p^{2\theta v}.
\ena
Thus
\bea\label{h(p,k) estimate}
|h_{\pi}(k,p)| &\ll_{m,r}& k^{r^2 m^2 + 2m - 2}
\biggl\{\sum_{v=1}^{k-1}v^{m^2}p^{2\theta v}
(k-v)^{m^2}p^{2\theta (k-v)}+k^{m^2}p^{2\theta k}\biggr\}\nonumber\\
&\ll_{m,r}& k^{C_1}p^{2\theta k}.
\eea
where $C_1=C_1(m,r)=r^2 m^2 + 2m^2+2m - 1$.

Now we estimate the contributions from
terms with $k\geq 2$  to \eqref{H_p}. For $k\geq m^{*2}+2$, we
set
\bna
\sigma_0=1-\frac{1}{m^{*2}+1}+\varepsilon.
\ena
By \eqref{h(p,k) estimate}, for $\sigma\geq \sigma_0$,
\bea\label{k large}
\sum_{k\geq m^{*2}+2}\frac{|h_{\pi}(k,p)|}{p^{k\sigma}}
\ll \sum_{k\geq m^{*2}+2}k^{C_1}p^{-k(\sigma-2\theta)}
\ll_{\pi,\varepsilon} p^{-1-\frac{1}{m^{*2}+1}}.
\eea
Here $\theta=1/2-1/(m^{*2}+1)$ is defined in \eqref{max theta}.

By \eqref{eq:linear-cancellation}, \eqref{k large} and
the trivial inequality $1+a+b\leq (1+a)(1+b)$,
for $\mathrm{Re}(s)\geq \sigma_0$, the product over the unramified primes
in
\eqref{H_p} can be bounded by
 \bea\label{H-pi-sigma}
\prod_{p\nmid q_{\pi}}\Big(1+
       \sum_{k\ge 1}\frac{|h_{\pi}(k,p)|}{p^{k\sigma}}
       \Big)&\ll& \prod_{p\nmid q_{\pi}}\Big(1+
       \sum_{k=2}^{m^{*2}+1}\frac{|h_{\pi}(k,p)|}{p^{k\sigma}}
       \Big)\prod_{p\nmid q_{\pi}}\Big(1+
       \sum_{k\geq m^{*2}+2}\frac{|h_{\pi}(k,p)|}{p^{k\sigma}}
       \Big)\nonumber\\
&\ll&\prod_{k=2}^{m^{*2}+1}\prod_{p\nmid q_{\pi}}\Big(1+\frac{|h_{\pi}(k,p)|}{p^{k\sigma}}
\Big).
\eea

For the remaining range $2\leq k\leq m^{*2}+1$, we first reduce the convolution term in \eqref{eq:h-local-majorant} to a sum involving a single Rankin--Selberg coefficient. Indeed, by the bound $B(v)\ll_r v^{m^2}p^{2\theta v}$ established above, we have
\begin{align*}
\sum_{v=1}^{k-1}B(v)B(k-v)
&\ll_{m,r}\sum_{v=1}^{k-1}B(v)p^{2\theta(k-v)}\\
&\ll_{m,r}\sum_{j=1}^r\sum_{\ell=1}^{k-1}
\lambda_{\pi_j\times\widetilde\pi_j}(p^\ell)p^{2\theta(k-\ell)}.
\end{align*}
 Combining this with \eqref{eq:h-local-majorant}, we obtain
\begin{align*}
|h_\pi(k,p)|\ll_{m,r}
\sum_{j=1}^r\sum_{\ell=1}^k
\lambda_{\pi_j\times\widetilde\pi_j}(p^\ell)p^{2\theta(k-\ell)}.
\end{align*}

We now separate the terms with $\ell=1$ from those with $2\leq\ell\leq k$. Since $2\theta=1-2/(m^{*2}+1)$, for $k\geq2$ and $\sigma\geq\sigma_0$ we have
\begin{align}\label{eq:linear-local-exponent}
k\sigma-2\theta(k-1)
&\geq 1+\frac{k-2}{m^{*2}+1}+k\varepsilon
\geq 1+\varepsilon.
\end{align}
Consequently,
\begin{align*}
\sum_{k=2}^{m^{*2}+1}
\frac{\lambda_{\pi_j\times\widetilde\pi_j}(p)p^{2\theta(k-1)}}{p^{k\sigma}}
\ll_{m,r}
\frac{\lambda_{\pi_j\times\widetilde\pi_j}(p)}{p^{1+\varepsilon}}.
\end{align*}
For $2\leq\ell\leq k$, since $\sigma\geq\sigma_0>2\theta$, we similarly obtain
\begin{align*}
k\sigma-2\theta(k-\ell)
=\ell\sigma+(k-\ell)(\sigma-2\theta)\geq\ell\left(1-\frac{1}{m^{*2}+1}+\varepsilon\right).
\end{align*}
It follows that
\begin{align*}
\sum_{k=2}^{m^{*2}+1}\sum_{\ell=2}^k
\frac{\lambda_{\pi_j\times\widetilde\pi_j}(p^\ell)p^{2\theta(k-\ell)}}{p^{k\sigma}}\ll_{m,r}
\sum_{\ell=2}^{\infty}
\frac{\lambda_{\pi_j\times\widetilde\pi_j}(p^\ell)}
{p^{\ell(1-1/(m^{*2}+1)+\varepsilon)}}.
\end{align*}
Combining these estimates and using the nonnegativity of the Rankin--Selberg coefficients, we deduce that
\begin{align*}
\prod_{p\nmid q_\pi}\left(1+\sum_{k=2}^{m^{*2}+1}
\frac{|h_\pi(k,p)|}{p^{k\sigma}}\right)
\leq&
\prod_{j=1}^r\left\{
\prod_p\left(1+
\frac{\lambda_{\pi_j\times\widetilde\pi_j}(p)}{p^{1+\varepsilon}}
\right)\right\}^{O_{m,r}(1)}\\ &\times
\prod_{j=1}^r\left\{
\prod_p\left(1+\sum_{\ell=2}^{\infty}
\frac{\lambda_{\pi_j\times\widetilde\pi_j}(p^\ell)}
{p^{\ell(1-1/(m^{*2}+1)+\varepsilon)}}
\right)\right\}^{O_{m,r}(1)}.
\end{align*}
Here we have used the elementary inequality
$1+C\sum_k x_k\leq\prod_k(1+x_k)^C$ for $x_k\geq0$ and $C\geq1$.
The first Euler product is bounded by $L(1+\varepsilon,\pi_j\times\widetilde\pi_j)$,
while the second is controlled by Theorem~\ref{Jiang}, since $m_j\leq m^*$ implies
\[
1-\frac{1}{m^{*2}+1}+\varepsilon
\geq 1-\frac{1}{m_j^2+1}+\varepsilon.
\]
Therefore,
\begin{align*}
\prod_{p\nmid q_\pi}\left(1+\sum_{k=2}^{m^{*2}+1}
\frac{|h_\pi(k,p)|}{p^{k\sigma}}\right)\ll_{m,r,\varepsilon}
\prod_{j=1}^r
L(1+\varepsilon,\pi_j\times\widetilde\pi_j)^{O_{m,r}(1)}
\prod_{j=1}^r\mathcal C(\pi_j)^{O_{m,r}(\varepsilon)}.
\end{align*}
Together with \eqref{k large} and \eqref{H-pi-sigma}, this proves absolute
and locally uniform convergence of the unramified correction product.
The complete factor is
\[
H_\pi(s)=\prod_{p\mid q_\pi}H_{\pi,p}(s)
 \prod_{p\nmid q_\pi}\left(1+\sum_{k\ge2}\frac{h_\pi(k,p)}{p^{ks}}\right).
\]
The finitely many factors $H_{\pi,p}(s)=D_p(s)L_p(s,\pi\times\widetilde\pi)^{-1}$
for $p\mid q_\pi$ are holomorphic and bounded in the stated closed
half-planes, since these lie in $\operatorname{Re}s>2\theta$.
Consequently $H_\pi(s)\ll_{\pi,\varepsilon}1$ there. All representations
are fixed; no uniformity in a varying conductor is claimed.

For the second factorization, write
\[
 \mathcal L_p(X)=\prod_{i,j=1}^r
 L_p(s,\pi_i\times\widetilde\pi_j)\bigm|_{p^{-s}=X}.
\]
For $p\nmid q_\pi$, by \eqref{eq:linear-cancellation}, the coefficient of $X$ in
$(1-|\lambda_\pi(p)|^2X)\mathcal L_p(X)$ is zero.  Its coefficients of degree at least
two satisfy the same majorants as \eqref{eq:h-local-majorant}: the terms of
degree at least two in $(1-|\lambda_\pi(p)|^2X)\mathcal L_p(X)$ are covered by Theorem~\ref{Jiang}, and
the remaining product with $|\lambda_\pi(p)|^2X$ is covered by
\eqref{eq:linear-local-exponent}.  Hence
\[
 H_{\pi}^{\sharp}(s)=\prod_p
 \left(1-\frac{|\lambda_{\pi}(p)|^2}{p^s}\right)
 \prod_{i,j=1}^rL_p(s,\pi_i\times\widetilde\pi_j)
\]
is also holomorphic and bounded in the asserted closed half-planes.
The unconditional argument uses Theorem~\ref{Jiang} to bound the higher
prime-power coefficients on average for each cuspidal constituent.
Under GRC, cancellation of the linear terms at unramified primes gives
\(H_{\pi,p}(s)=1+O_{m,\varepsilon}(p^{-2\sigma})\) and
\(H_{\pi,p}^{\sharp}(s)=1+O_{m,\varepsilon}(p^{-2\sigma})\) uniformly
for \(\sigma=\Re s\geq 1/2+\varepsilon\). The finitely many ramified
factors are holomorphic there by the tempered local bounds. The
conditional assertion follows by summing \(p^{-2\sigma}\) over primes.
\end{proof}

The lower coefficients are handled in the same way, one derivative at a
time.  For \(\mathrm{Re}(z)>1\) and \(s\) in a neighborhood of \(1\), normal
convergence permits termwise differentiation.  Hence, for each fixed
\(0\leq k\leq r-1\),
\begin{align}\label{eq:derivative-modulus-series}
 \sum_{q\geq1}\frac{V_q^{(k)}(1)}{q^z}
 &=
 \left.
 \frac{\partial^k}{\partial s^k}
 \prod_p
 \left(
  1+D_p(s)^{-1}\sum_{\ell\geq1}p^{-\ell z}
 \right)
 \right|_{s=1}\notag\\
 &=
 \zeta(z)B_{\pi,k}(z),
\end{align}
where
\begin{align*}
 B_{\pi,k}(z)
 :=
 \left.
 \frac{\partial^k}{\partial s^k}
 \prod_p
 \left\{
  1+\bigl(D_p(s)^{-1}-1\bigr)p^{-z}
 \right\}
 \right|_{s=1}.
\end{align*}
It then follows \(B_{\pi,k}(z)\) is holomorphic for
\(\mathrm{Re}(z)>2\theta\)  and is bounded in every closed half-plane contained
in this region.
Applying Mellin inversion to
\eqref{eq:derivative-modulus-series} and shifting the \(z\)-contour across
the simple pole of \(\zeta(z)\), we obtain
\begin{align}\label{eq:average-Vq-derivative}
 \sum_q\Phi(q/Q)V_q^{(k)}(1)
 =
 QB_{\pi,k}(1)+
 O_{\pi,\Phi,k,\varepsilon}
 \left(Q^{2\theta+\varepsilon}\right).
\end{align}
Leibniz's rule gives
 \[ u_{q,\ell}
 =
 \frac{1}{\ell!}
 \sum_{k=0}^{\ell}
 \binom{\ell}{k}
 A_\pi^{(\ell-k)}(1)V_q^{(k)}(1).
 \]
Consequently, if
\begin{align}\label{eq:beta-Pi-nu}
 \beta_{\pi,\ell}
 :=
 \frac{1}{\ell!}
 \sum_{k=0}^{\ell}
 \binom{\ell}{k}
 A_\pi^{(\ell-k)}(1)B_{\pi,k}(1),
\end{align}
then \eqref{eq:average-Vq-derivative} implies
\begin{align}\label{eq:average-uq-nu}
 \sum_q\Phi(q/Q)u_{q,\ell}
 =
 Q\beta_{\pi,\ell}+
 O_{\pi,\Phi,\Psi,\ell,\varepsilon}
 \left(Q^{2\theta+\varepsilon}\right).
\end{align}
Define
\begin{align}\label{eq:fixed-diagonal-polynomial}
 P_{r-1}(Y)
 :=
 \sum_{\ell=0}^{r-1}
 \frac{\beta_{\pi,\ell}}{(r-1-\ell)!}
 Y^{r-1-\ell}.
\end{align}
Since \(B_{\pi,0}(1)>0\) and \(A_\pi(1)\ne0\), the polynomial
\(P_{r-1}\) has degree exactly \(r-1\).
Equations \eqref{eq:diagonal-polynomial} and
\eqref{eq:average-uq-nu} then give
 \begin{align}\label{eq:averaged-diagonal-residue}
 \sum_q\Phi(q/Q)XP_{q,r-1}(\log X)
 &=
 QXP_{r-1}(\log X)\notag\\
 &\quad+
 O_{\pi,\Phi,\Psi,\varepsilon}
 \left(
  Q^{2\theta+\varepsilon}
  X(1+\log X)^{r-1}
 \right).
\end{align}

It remains to estimate the integral on the shifted line in
\eqref{eq:4.7}.  Put
$\sigma_0=1-1/(m^{*2}+1)+\varepsilon$.
Then
\begin{align*}
 \sigma_0-2\theta
 \geq
 \frac{1}{m^{*2}+1}+\varepsilon>0.
\end{align*}
Uniformly for \(\mathrm{Re}(s)=\sigma_0\), the local coefficient bounds imply
\begin{align*}
 D_p(s)-1
 &\ll_\pi
 \sum_{\ell\geq1}
 (\ell+1)^{O_\pi(1)}
 p^{-\ell(\sigma_0-2\theta)}\\
 &\ll_{\pi,\varepsilon}
 p^{-1/(m^{*2}+1)-\varepsilon}.
\end{align*}
Consequently, for all sufficiently large \(p\),
\begin{align}\label{eq:local-R-inverse-bound}
 D_p(s)^{-1}
 =
 1+
 O_{\pi,\varepsilon}
 \left(
  p^{-1/(m^{*2}+1)-\varepsilon}
 \right).
\end{align}
There are only finitely many remaining primes.  After altering the real part
of the contour slightly, if necessary, and then relabelling
\(\varepsilon\), we may assume that none of the corresponding local factors
vanishes on the chosen line; their inverses are then uniformly bounded
there.  It follows from \eqref{eq:local-R-inverse-bound} and the standard
divisor bound that
\begin{align*}
 \prod_{p\mid q}
 \left|D_p(s)^{-1}\right|
 \ll_{\pi,\varepsilon}q^\varepsilon
\end{align*}
uniformly for \(\mathrm{Re}(s)=\sigma_0\).

We divide the vertical integral into the range \(|\mbox{Im}(s)|\leq1\) and the
dyadic ranges
\[
 T\leq|\mbox{Im}(s)|\leq2T.
\]
The Rankin--Selberg factors have at most polynomial growth in vertical
strips, while \(\widetilde{\Psi^2}(s)\) decays faster than any power of
\(|\mbox{Im}(s)|\).  Lemma~\ref{corollary:Jiang} gives the required uniform
bound for \(H_\pi(s)\) on this line.  Hence
\[
 \int_{(\sigma_0)}
 \left|
  \widetilde{\Psi^2}(s)H_\pi(s)
  \prod_{h,j=1}^r
  L(s,\pi_h\times\widetilde\pi_j)
  \prod_{p\mid q}D_p(s)^{-1}
 \right|
 |\dif s|
 \ll_{\pi,\Psi,\varepsilon}q^\varepsilon.
\]
After multiplication by \(X^{\sigma_0}\) and summation over
\(q\asymp Q\), the shifted integral contributes
\begin{align}\label{eq:diagonal-shifted-integral}
 O_{\pi,\Phi,\Psi,\varepsilon}
 \left(
  QX^{1-1/(m^{*2}+1)+\varepsilon}
 \right).
\end{align}
In the range under consideration, the harmless factor \(Q^\varepsilon\)
may be absorbed into \(X^\varepsilon\), after changing \(\varepsilon\).

Combining \eqref{eq:averaged-diagonal-residue} and
\eqref{eq:diagonal-shifted-integral}, we obtain
\begin{align}\label{eq:final-diagonal-asymptotic}
\mathscr A_\pi(Q,X)
 &=
 QXP_{r-1}(\log X)\notag\\
 &\quad+
 O_{\pi,\Phi,\Psi,\varepsilon}
 \left(
  Q^{2\theta+\varepsilon}
  X^{1+\varepsilon}
  +
  QX^{1-1/(m^{*2}+1)+\varepsilon}
 \right).
\end{align}
Since \(\theta<1/2\), choosing \(\varepsilon>0\) sufficiently small shows
that the error term in \eqref{eq:final-diagonal-asymptotic} is \(o(QX)\)
uniformly in the range of Proposition~\ref{finald}, where \(X\) and \(Q\) are polynomially
related.

Finally, \eqref{eq:q-leading-coefficient},
\eqref{eq:average-Vq-zero}, and \eqref{eq:fixed-diagonal-polynomial} show
that the leading coefficient of \(P_{r-1}\) is
\[
 \frac{B_{\pi,0}(1)H_\pi(1)}{(r-1)!}
 \prod_{h=1}^r
 \operatorname*{Res}_{s=1}
 L(s,\pi_h\times\widetilde\pi_h)
 \prod_{\substack{h,j=1\\h\ne j}}^r
 L(1,\pi_h\times\widetilde\pi_j),
\]
which is precisely the constant in \eqref{eq:leading-constant}.
Moreover, \eqref{eq:beta-Pi-nu} and
\eqref{eq:fixed-diagonal-polynomial} give every lower-order coefficient
explicitly.  This completes the proof of Proposition~\ref{finald}.

\section{Proof of Proposition~\ref{finalg1}}
\label{sec:proof-finalg1}

We estimate \( \mathscr B_{1,\pi}(Q,X)\) by expressing it as a triple
Mellin integral and then shifting the Mellin contours.  The continuation
needed for these shifts is provided by the following factorization from
Lemma~\ref{analyticcontinuationofTk}:
\[
 T(s_1,s_2;z)
 =\zeta(1-z)
 \frac{\displaystyle
  L(s_1,\widetilde\pi)L(s_2,\pi)
  \prod_{h,j}
  L(s_1+s_2-z,\pi_h\times\widetilde\pi_j)}
 {\displaystyle
  \prod_{h,j}
  L(s_1+s_2,\pi_h\times\widetilde\pi_j)}
 \mathcal P(s_1,s_2,z)
\]
throughout the region stated there.  On applying Mellin inversion to
\(\mathcal W\), we therefore obtain, for any \(A>0\),
\[
  \mathscr B_{1,\pi}(Q,X)
 =\frac{1}{(2\pi i)^3}
 \int_{(-A)}
 \int_{(1+\varepsilon)}
 \int_{(1+\varepsilon)}
 \widetilde{\mathcal W}_3(s_1,s_2;z)
 X^{s_1+s_2}
 \left(\frac QX\right)^z
 T(s_1,s_2;z)
 \dif s_1\dif s_2\dif z.
\]

Then we move the \(s_1\)  and
\(s_2\) contours into the region supplied by
Lemma~\ref{analyticcontinuationofTk}, while choosing them so that the
Rankin--Selberg factors depending on \(s_1+s_2\) remain in their regions of
absolute convergence.  The decay of
\(\widetilde{\mathcal W}_3\), together with mean-value estimates for
\(L(s,\pi)\), will then control the resulting integrals. Precisely, for fixed $z$ on $\mathrm{Re}(z)=-A$, we move both $s$ contours to the contour
$\Gamma$ consisting of the five line segments
\[
 (1+\varepsilon)-i\infty,\quad
 (1+\varepsilon)-iY,\quad
 (5/6+\varepsilon)-iY,\quad
 (5/6+\varepsilon)+iY,\quad
 (1+\varepsilon)+iY,\quad
 (1+\varepsilon)+i\infty.
\]
 Let $V$ be the central vertical segment of $\Gamma$, let $H^\pm$ be its two
horizontal segments, and let $V_\infty^\pm$ be the two vertical tails.

The standard factors $L(s,\pi)$ are entire, and
Lemma~\ref{analyticcontinuationofTk} justifies the deformation.  If
$s_1,s_2\in\Gamma$ and $\mathrm{Re}(z)=-A$, then
\[
 \mathrm{Re}(s_1+s_2-z)\geq\frac53+A+2\varepsilon,
 \qquad
 \mathrm{Re}(s_1+s_2)\geq\frac53+2\varepsilon.
\]
The Rankin--Selberg $L$-functions in both the numerator and denominator are
therefore represented by absolutely convergent Euler products.  The same is
true for $\zeta(1-z)$, while $\mathcal P(s_1,s_2,z)$ is bounded by
Lemma~\ref{analyticcontinuationofTk}.  Hence the only factors requiring mean
value estimates are $L(s_1,\widetilde\pi)$ and $L(s_2,\pi)$. Write $s_k=\sigma_k+it_k$ and $z=-A+i\tau$. Taking an integer $P\geq 1+A$ in
Proposition~\ref{MultiMellinBound}, the absolute value of the integrand in
the two $s$-variables is bounded by
\begin{align*}
 \frac{X^{\sigma_1+\sigma_2+A}
 |L(s_1,\widetilde\pi)L(s_2,\pi)|}
 {(1+|\tau|^{1+A})(1+|t_1+t_2|^{1+A})
  (1+|t_1-t_2|^{1+A})}.
\end{align*}

We require a convolution estimate along the entire contour \(\Gamma\).
The vertical components \(V\cup V_\infty^+\cup V_\infty^-\) together
cover the full range \(t_2\in\mathbb R\).  For fixed \(t_1\), the change
of variables \(u=t_1+t_2\) gives
\begin{align*}
 &\int_{-\infty}^{\infty}
 \frac{\dif t_2}
 {(1+|t_1+t_2|^{1+A})(1+|t_1-t_2|^{1+A})}\\
 &\qquad=
 \int_{-\infty}^{\infty}
 \frac{\dif u}
 {(1+|u|^{1+A})(1+|2t_1-u|^{1+A})}
 \ll_A
 \frac{1}{1+|t_1|^{1+A}}.
\end{align*}
Indeed, the real line is covered by the two regions
\[
 |u|\geq|t_1|
 \qquad\text{and}\qquad
 |2t_1-u|\geq|t_1|.
\]
On each region, one factor is
\(O_A((1+|t_1|^{1+A})^{-1})\), while the other is integrable.

It remains to account for the horizontal components of \(\Gamma\).
On \(H^\pm\), we have \(t_2=\pm Y\), and the total length of these
segments is \(O(1)\).  Since
\[
 \max\{|t_1+Y|,|t_1-Y|\}\geq |t_1|,
\]
we have
\begin{align*}
 \int_{H^+\cup H^-}
 \frac{|\dif s_2|}
 {(1+|t_1+t_2|^{1+A})(1+|t_1-t_2|^{1+A})}
 \ll_A
 \frac{1}{1+|t_1|^{1+A}}.
\end{align*}
Combining the vertical and horizontal contributions yields
\begin{align}\label{eq:kernel-convolution}
 \int_\Gamma
 \frac{|\dif s_2|}
 {(1+|t_1+t_2|^{1+A})(1+|t_1-t_2|^{1+A})}
 \ll_A
 \frac{1}{1+|t_1|^{1+A}}.
\end{align}
The analogous estimate with \(s_1\) and \(s_2\) interchanged also holds.

Since no poles are crossed, shifting both Mellin contours produces no
residue terms, and the resulting double integral is taken over
\(\Gamma\times\Gamma\).  On this product contour, we use
\[
 X^{\sigma_1+\sigma_2}
 |L(s_1,\widetilde\pi)L(s_2,\pi)|
 \leq
 \frac12X^{2\sigma_1}|L(s_1,\widetilde\pi)|^2
 +\frac12X^{2\sigma_2}|L(s_2,\pi)|^2.
\]
For the first term, we integrate with respect to \(s_2\) using
\eqref{eq:kernel-convolution}; for the second, we integrate with respect
to \(s_1\).  The two resulting integrals are estimated identically, since
the same bounds apply to \(\pi\) and \(\widetilde\pi\).  We are therefore
reduced to
\begin{align}\label{eq:G1-one-variable-reduction}
 \frac{1}{1+|\tau|^{1+A}}
 \int_\Gamma
 \frac{X^{2\mathrm{Re}(s)+A}|L(s,\pi)|^2}
      {1+|\mbox{Im}(s)|^{1+A}}
 |\dif s|.
\end{align}
Here
\[
 L(s,\pi)
 =\prod_{j=1}^rL(s,\pi_j)
\]
is an automorphic \(L\)-function of degree \(m\).

We estimate separately the vertical tails, the horizontal segments, and
the central segment of \(\Gamma\).  On \(V_\infty^\pm\), we have
\(\mathrm{Re}(s)=1+\varepsilon\), so absolute convergence gives
\[
 L(s,\pi)\ll_{\pi,\varepsilon}1.
\]
It follows that
\begin{align}\label{eq:G1-tail}
 \int_{V_\infty^+\cup V_\infty^-}
 \frac{X^{2\mathrm{Re}(s)+A}|L(s,\pi)|^2}
      {1+|\mbox{Im}(s)|^{1+A}}
 |\dif s|
 \ll_{\pi,\varepsilon}
 X^{2+A+\varepsilon}Y^{-A}.
\end{align}

On the horizontal segments
\[
 H^\pm
 =\left\{
  \sigma\pm iY:
  \frac56+\varepsilon\leq\sigma\leq1+\varepsilon
 \right\},
\]
the convexity bound of $L$-function gives
\begin{align}\label{eq:G1-horizontal}
 \int_{H^+\cup H^-}
 \frac{X^{2\mathrm{Re}(s)+A}|L(s,\pi)|^2}
      {1+|\mbox{Im}(s)|^{1+A}}
 |\dif s|
 \ll_{\pi,\varepsilon}{}&
 X^{2+A+\varepsilon}Y^{-1-A+\varepsilon}
\\
 &+
 X^{5/3+A+\varepsilon}
 Y^{m/6-1-A+\varepsilon}.\nonumber
\end{align}

It remains to estimate the central segment
\[
 V=\left\{
  \frac56+\varepsilon+it:|t|\leq Y
 \right\}.
\]
Lemma~\ref{meanLbound} gives, for \(U\geq1\),
\[
 \int_U^{2U}
 \left|
  L\left(\frac56+\varepsilon+it,\pi\right)
 \right|^2
 \dif t
 \ll_{\pi,\varepsilon}
 U^{m/6+2/3+\varepsilon}.
\]
After inserting the factor
\((1+|t|^{1+A})^{-1}\), the contribution of the dyadic range
\(U\leq|t|\leq2U\) is
\[
 \ll_{\pi,\varepsilon}
 X^{5/3+A+\varepsilon}
 U^{m/6-1/3-A+\varepsilon}.
\]
The range \(|t|\leq1\) contributes
\(O_{\pi,\varepsilon}(X^{5/3+A+\varepsilon})\).  Summing the dyadic
ranges \(1\leq U\leq Y\), we conclude that
\begin{align}\label{eq:G1-central}
 \int_V
\frac{X^{2\mathrm{Re}(s)+A}|L(s,\pi)|^2}
      {1+|\mbox{Im}(s)|^{1+A}}
 |\dif s|
 \ll_{\pi,\varepsilon}
 X^{5/3+A+\varepsilon}
  Y^{\max\{m/6-1/3-A,\,0\}+\varepsilon}
 .
\end{align}

The factor \((1+|\tau|^{1+A})^{-1}\) in
\eqref{eq:G1-one-variable-reduction} is integrable.  Therefore, integrating
over \(z=-A+i\tau\) and combining
\eqref{eq:G1-tail}--\eqref{eq:G1-central}, we obtain
\begin{align*}
  \mathscr B_{1,\pi}(Q,X)
\ll_{\pi,\Phi,\Psi,\varepsilon}
 Q^{-A}\Big(
  X^{2+A+\varepsilon}Y^{-A}
  +X^{5/3+A+\varepsilon}
    Y^{\max\{m/6-1/3-A,\,0\}+\varepsilon}
 \Big).
\end{align*}
Here the two terms in \eqref{eq:G1-horizontal} have been absorbed into
the corresponding tail and central-segment estimates.

Suppose first that \(m>2\).  Take
\[
 A=\frac{m-2}{6},
 \qquad
 Y=X^{2/(m-2)}.
\]
The exponent of \(Y\) arising from the central segment is then zero, and
the resulting logarithmic loss is absorbed into \(X^\varepsilon\).
Moreover,
\[
 X^{2+A}Y^{-A}=X^{5/3+A}.
\]
Hence
\[
  \mathscr B_{1,\pi}(Q,X)
 \ll_{\pi,\Phi,\Psi,\varepsilon}
 Q^{-(m-2)/6}
 X^{5/3+(m-2)/6+\varepsilon}.
\]
This is \(o(QX)\) throughout the asserted range
\[
 X\leq Q^{1+2/(m+2)-\varepsilon}.
\]

If $m=2$, the assumptions imply $r=1$ and $m_1=2$.
The Kim--Sarnak bound \cite[Appendix~2]{ks} permits the local exponent
$\vartheta=7/64<1/6$ in Lemma~\ref{analyticcontinuationofTk}.
We may therefore place the central segment on
\(\mathrm{Re}(s)=1/2+\varepsilon\).  Applying Lemma~\ref{meanLbound} with parameter \(\varepsilon/2\) gives
\[
 \int_U^{2U}
 |L(1/2+\varepsilon+it,\widetilde\pi)|^2
 \dif t
 \ll_{\pi,\varepsilon}U^{1+\varepsilon}.
\]
Repeating the preceding argument gives, for every \(A>0\),
\[
  \mathscr B_{1,\pi}(Q,X)
 \ll_{\pi,\Phi,\Psi,\varepsilon}
 Q^{-A}\left(
  X^{2+A+\varepsilon}Y^{-A}
  +X^{1+A+\varepsilon}
 \right).
\]
Letting \(Y\to\infty\) and then choosing \(A>0\) sufficiently small shows
that
\[
  \mathscr B_{1,\pi}(Q,X)=o(QX)
\]
uniformly for
\[
 X\leq Q^{2-\varepsilon}.
\]

Assume finally GRH and GRC as in Theorem~\ref{mainresult2} .
The latter permits $\vartheta=0$ in Lemma~\ref{analyticcontinuationofTk},
so for every $m$ the central segment may be placed on
\(\mathrm{Re}(s)=1/2+\varepsilon\).  GRH gives the Lindel\"of bound
\[
 L(1/2+\varepsilon+it,\pi)
 \ll_{\pi,\varepsilon}(2+|t|)^\varepsilon.
\]
The same decomposition therefore yields
\[
  \mathscr B_{1,\pi}(Q,X)
 \ll_{\pi,\Phi,\Psi,\varepsilon}
 Q^{-A}\left(
  X^{2+A+\varepsilon}Y^{-A}
  +X^{1+A+\varepsilon}
 \right).
\]
Letting \(Y\to\infty\) and taking \(A>0\) sufficiently small proves
\[
  \mathscr B_{1,\pi}(Q,X)=o(QX)
\]
uniformly for \(X\leq Q^{2-\varepsilon}\), as required.

\section{Analytic tools}

The decomposition in Section~3 reduces the main theorems to
Propositions~\ref{finald}, \ref{finalg1}, and~\ref{finalg2}.
Section~\ref{sec:proof-finald} treats the diagonal term by factoring the
square-coefficient Dirichlet series in terms of Rankin--Selberg \(L\)-functions.
In the off-diagonal terms, the variables \(n_1,n_2,d,g\) remain coupled
through the weight
\[
 \mathcal W\left(\frac{n_1}{X},\frac{n_2}{X};\frac{X}{dgQ}\right),
\]
Multivariable Mellin transforms separate the arithmetic variables in the
weight and provide the decay required for the contour shifts in
Sections~\ref{sec:proof-finalg1} and~\ref{sec:proof-finalg2}.
The principal-character contribution leads to the multiple Dirichlet series
\(T(s_1,s_2;z)\), whose meromorphic continuation and factorization identify
the automorphic \(L\)-functions governing \(\mathscr B_{1,\pi}(Q,X)\).
This factorization, used in Section~\ref{sec:proof-finalg1}, is stated in
Lemma~\ref{analyticcontinuationofTk} and proved in
Section~\ref{sec:proof-analyticcontinuationofTk}.

\subsection{Multivariable Mellin transforms}

Recall that the weight introduced in the off-diagonal decomposition is
\[
 \mathcal W(x,y;u)=\Psi(x)\Psi(y)\Phi(u|x-y|),
 \qquad x,y,u>0.
\]
This is the same notation and normalization as in
\cite[Section~3]{RS}.  We restate the transforms here so that the analytic
input used later is explicit.  For fixed $u>0$ define
\[
 \widetilde{\mathcal W}_2(s_1,s_2;u)
 =\int_0^\infty\int_0^\infty
  \mathcal W(x,y;u)x^{s_1-1}y^{s_2-1}\dif x\dif y,
\]
and define the full Mellin transform by
\[
 \widetilde{\mathcal W}_3(s_1,s_2;z)
 =\int_0^\infty\int_0^\infty\int_0^\infty
  \mathcal W(x,y;u)x^{s_1-1}y^{s_2-1}u^{z-1}
  \dif u\dif x\dif y.
\]
The corresponding inversion formulas are
\begin{align*}
 \mathcal W(x,y;u)
 &=\frac{1}{(2\pi i)^2}
   \int_{(\sigma_1)}\int_{(\sigma_2)}
   \widetilde{\mathcal W}_2(s_1,s_2;u)
   x^{-s_1}y^{-s_2}\dif s_2\dif s_1,\\
 \mathcal W(x,y;u)
 &=\frac{1}{(2\pi i)^3}
   \int_{(\sigma_1)}\int_{(\sigma_2)}\int_{(\sigma_3)}
   \widetilde{\mathcal W}_3(s_1,s_2;z)
   x^{-s_1}y^{-s_2}u^{-z}\dif z\dif s_2\dif s_1.
\end{align*}
Here $\sigma_1$ and $\sigma_2$ may be chosen freely, while
$\sigma_3<1$; the first two integrations are interpreted in the same
principal-value sense as in \cite[Section~3]{RS}.

We record only the two estimates from \cite[Section~3]{RodgersSoundararajan2018}
that will be used later.
\begin{proposition}\label{MultiMellinBound2}
Fix $A>0$.  Uniformly for $-A\leq\mathrm{Re}(s_1),\mathrm{Re}(s_2)\leq A$ and every integer
$P\geq1$,
\[
 \widetilde{\mathcal W}_2(s_1,s_2;u)
 \ll_{A,P}\frac{u^{P-1}}
 {1+\max\{|\Im(s_1)|,|\Im(s_2)|\}^{P}}.
\]
Moreover, there is a constant $K>0$, depending only on the supports of
$\Phi$ and $\Psi$, such that
$\widetilde{\mathcal W}_2(s_1,s_2;u)=0$ unless $u\geq K$.
\end{proposition}

\begin{proposition}\label{MultiMellinBound}
Fix $\varepsilon,A>0$.  Uniformly for
$-A\leq\mathrm{Re}(z)\leq1-\varepsilon$ and
$-A\leq\mathrm{Re}(s_1),\mathrm{Re}(s_2)\leq A$, one has, for every integer $P\geq1$,
\[
 \widetilde{\mathcal W}_3(s_1,s_2;z)
 \ll_{\varepsilon,A,P}
 \frac{1}{1+|\Im(z)|^P}
 \frac{1}{1+|\Im(s_1)+\Im(s_2)|^P}
 \frac{1}{1+|\Im(s_1)-\Im(s_2)|^{1-\mathrm{Re}(z)}}.
\]
\end{proposition}

These statements are quoted from \cite[Propositions~3.3 and 3.4]{RS} and
will be used with the displayed normalization.

\subsection{Meromorphic continuation of \(T\)}

The principal-character contribution
\( \mathscr B_{1,\pi}(Q,X)\), after Mellin inversion, leads naturally to
the multiple Dirichlet series
 \[T(s_1,s_2;z)
 :=
 \sum_{\substack{n_1=dN_1,\ n_2=dN_2\\(N_1,N_2)=1}}
 \sum_{g\geq1}
 \sum_{\substack{c\mid d\\(N_1N_2,cg)=1}}
 \frac{\overline{\lambda_\pi(n_1)}}{n_1^{s_1}}
 \frac{\lambda_\pi(n_2)}{n_2^{s_2}}
 \frac{\mu(c)}{\phi(cg)}
 (dg)^z.
\]
The series initially converges absolutely for $\operatorname{Re}s_1,
\operatorname{Re}s_2>1$ and $\operatorname{Re}z<0$.
We have interchanged the two coefficient variables, using the symmetry
of $\mathcal W$, so that the first coefficient is conjugated.
To shift the Mellin contours beyond the region of absolute convergence, we
require both a continuation of \(T(s_1,s_2;z)\) and an explicit
factorization in terms of automorphic \(L\)-functions.  The following lemma
provides this input.  Since its proof requires a lengthy prime-by-prime
analysis of the local factors of \(T\), we defer it to
Section~\ref{sec:proof-analyticcontinuationofTk}.

\begin{lem}\label{analyticcontinuationofTk}
    Let $\varepsilon>0$ be sufficiently small.  Unconditionally, the sum \(T(s_1,s_2;z)\) admits a meromorphic continuation to
$\mathrm{Re}(s_1),\mathrm{Re}(s_2)\geq\frac13+\theta+\varepsilon$ and $\mathrm{Re}(z)\leq2\varepsilon$, and in this region
   \begin{align*}
      T(s_1,s_2;z)=\zeta(1-z)
      \frac{\displaystyle
       \prod_hL(s_1,\widetilde\pi_h)
       \prod_jL(s_2,\pi_j)
       \prod_{h,j}L(s_1+s_2-z,\pi_h\times\widetilde\pi_j)}
      {\displaystyle\prod_{h,j}L(s_1+s_2,\pi_h\times\widetilde\pi_j)}
      \mathcal P(s_1,s_2,z),
   \end{align*}
   where $\mathcal P$ is a bounded, absolutely convergent Euler product.
   Under GRC the same assertion holds in the larger region
   \begin{align*}
    \{(s_1,s_2,z)\in\mathbb{C}^3 \mid\mathrm{Re}(s_1), \mathrm{Re}(s_2)\geq\frac{1}{3}+\varepsilon, \mathrm{Re}(z)\leq2\varepsilon\}.
   \end{align*}
\end{lem}

\section{Proof of Proposition~\ref{finalg2}}
\label{sec:proof-finalg2}

We now consider the contribution of the nonprincipal characters.  A direct
treatment of the resulting Mellin integrals would require moment estimates
for general automorphic \(L\)-functions that are not presently available in
the required strength.  To avoid this difficulty, we follow the strategy of
Rodgers and Soundararajan \cite{RS}.

\subsection{A hybrid mean-value estimate}
The estimate required for this argument is the following analogue of
\cite[Proposition~8]{RS}.

\begin{lem}\label{lem:hybridlargesieve}
For $R,S\geq1$,
\[
 \sum_{q\leq R}\ \sum_{\substack{\chi\bmod q\\\chi\ne\chi_0}}
 \int_{-S}^{S}|L(1/2+it,\pi\otimes\chi)|^2\dif t
 \ll_{\pi,\varepsilon}
 \bigl(R^2S+(RS)^{m/2}\bigr)(RS)^\varepsilon.
\]
Under GRH in Theorem~\ref{mainresult2}, this bound improves to $R^2S(RS)^\varepsilon$.
\end{lem}

\begin{proof}
First suppose that $\chi$ is primitive of conductor $q$.
Write $C_\chi(t)=\mathcal C(\Pi_\chi,t)$, so that
$C_\chi(t)\ll_\pi(q(2+|t|))^m$ by \eqref{eq:isobaric-L-degree&conductor}.
At $s_t=1/2+it$ the approximate functional equation
\cite[Theorem~5.3]{IK04} gives
\begin{align}\label{eq:twist-afe}
 L(s_t,\pi\times\chi)
 ={}&\sum_{n\ge1}\frac{\lambda_{\Pi_\chi}(n)}{n^{s_t}}\,
    V_{t,\chi}\!\left(\frac{n}{\sqrt{C_\chi(t)}}\right) \notag\\
 &+\varepsilon_{t,\chi}
   \sum_{n\ge1}\frac{\overline{\lambda_{\Pi_\chi}(n)}}{n^{1-s_t}}\,
    V^*_{t,\chi}\!\left(\frac{n}{\sqrt{C_\chi(t)}}\right),
 \qquad |\varepsilon_{t,\chi}|=1 .
\end{align}
With the normalization used in \eqref{eq:twist-afe},
\cite[Proposition~5.4]{IK04} implies that, for every \(B>0\) and every
integer \(k\geq0\),
\[
 y^k\left(
  \left|V_{t,\chi}^{(k)}(y)\right|
  +\left|(V^*_{t,\chi})^{(k)}(y)\right|
 \right)
 \ll_{\pi,B,k}(1+y)^{-B},
 \qquad y>0.
\]
Write $n=bv$ with $b\mid q_\pi^\infty$ and $(v,q_\pi)=1$, and put
$\vartheta=\max_j\theta_{m_j}<1/2$. The local factors at $p\nmid q_\pi$ and \eqref{theta} give
\[
 \lambda_{\Pi_\chi}(bv)=\lambda_{\Pi_\chi}(b)\lambda_\pi(v)\chi(v),
 \qquad |\lambda_{\Pi_\chi}(b)|\le d_m(b)b^\vartheta .
\]
Write $\|F\|_2^2=\sum_{q\le R}\sideset{}{^*}\sum_{\chi\bmod q}
\int_{-S}^S|F(\chi,t)|^2\dif t$.
For $Y\ge1$, put
\[
 \mathcal S_{\chi,t}(Y)=\sum_{\substack{v\le Y\\(v,q_\pi)=1}}
       \frac{\lambda_\pi(v)\chi(v)}{v^{s_t}}.
\]
Gallagher's hybrid large sieve \cite[Theorem~7.17]{IK04} and
Lemma~\ref{corollary:Jiang} give
\[
 \|\mathcal S_{\chi,t}(Y)\|_2^2
 \ll (R^2S+Y)\sum_{\substack{v\le Y\\(v,q_\pi)=1}}
       \frac{|\lambda_\pi(v)|^2}{v}
 \ll_{\pi,\eta}(R^2S+Y)Y^{2\eta}.
\]
Thus, for $x\ge1$, the triangle inequality gives
\begin{align}
 \left\|\sum_{n\le x}\frac{\lambda_{\Pi_\chi}(n)}{n^{s_t}}\right\|_2
 &=\left\|\sum_{\substack{b\le x\\ b\mid q_\pi^\infty}}
       \frac{\lambda_{\Pi_\chi}(b)}{b^{s_t}}\mathcal S_{\chi,t}(x/b)\right\|_2 \notag\\
 &\le \sum_{\substack{b\le x\\ b\mid q_\pi^\infty}}
       \frac{d_m(b)}{b^{1/2-\vartheta}}\|\mathcal S_{\chi,t}(x/b)\|_2 \notag\\
 &\ll_{\pi,\eta}(R^2S+x)^{1/2}x^\eta
        \sum_{b\mid q_\pi^\infty}\frac{d_m(b)}{b^{1/2-\vartheta}}.
 \label{eq:Gallagher-block}
\end{align}
The last sum equals $\prod_{p\mid q_\pi}(1-p^{\vartheta-1/2})^{-m}
\ll_\pi 1$, since $\vartheta<1/2$.
By rapid decay we may truncate \eqref{eq:twist-afe} at $N=(2RS)^{m/2+\eta}$
with error $O_{\pi,\eta,B}((2RS)^{-B})$.
Partial summation and \eqref{eq:Gallagher-block}, also applied to the
conjugate sum, give
\[
 \|L(s_t,\pi\times\chi)\|_2^2
 \ll_{\pi,\eta}(1+\log N)^2(R^2S+N)N^{2\eta}
 \ll_{\pi,\varepsilon}\bigl(R^2S+(RS)^{m/2}\bigr)(RS)^\varepsilon
\]
for $\eta>0$ sufficiently small in terms of $\varepsilon$.

For a general $\chi\bmod q$ induced by $\chi^*\bmod f$, Section~2.2
and \eqref{theta} give, uniformly in $t$,
\[
 |P_{\pi,\chi;q}(1/2+it)|
 \le\prod_{p\mid q}(1+p^{\vartheta-1/2})^m
 \ll_{\pi,\eta}q^\eta .
\]
Write $q=kf$. For each $k$ every primitive character of conductor
$f\le R/k$ induces a unique character modulo $kf$. Applying the
primitive estimate with $R/k$ in place of $R$ leaves the two sums
\[
 R^2 S\sum_{k\le R}k^{-2}
 \quad\hbox{and}\quad
 (RS)^{m/2}\sum_{k\le R}k^{-m/2}.
\]
Since $m\ge2$, these cost at most a logarithm, which is absorbed
in $(RS)^\varepsilon$ together with the factor from $P_{\pi,\chi;q}$.
This proves the asserted bound for the coefficientwise twists.

Under the stated GRH,
$L(1/2+it,\pi\times\chi^*)\ll_{\pi,\eta}(f(2+|t|))^\eta$
by the conductor-aspect Lindel\"of bound applied to each constituent.
Multiplying by $P_{\pi,\chi;q}$ and summing over $q\le R$ and
$\chi\bmod q$ gives $O_{\pi,\varepsilon}(R^2S(RS)^\varepsilon)$.
\end{proof}

\subsection{Dyadic decomposition}

After writing $n_1=dN_1$, $n_2=dN_2$, $d=cb$ and removing $(N_1,N_2)=1$
by M\"obius inversion, Mellin
inversion in the first two variables gives
\begin{align}\label{eq:G2-Mellin}
 \mathscr B_{2,\pi}(Q,X)
 =\sum_{\substack{c,b,e\ll X\\g\ll X/(cbQ)\\(e,cg)=1}}
 \frac{\mu(c)\mu(e)}{\phi(cg)}
 \sum_{\substack{\chi\bmod cg\\\chi\ne\chi_0}}
 \frac{1}{(2\pi i)^2}\int_{(\beta)}\int_{(\beta)}
 \frac{X^{s_1+s_2}}{(cbe)^{s_1+s_2}}\\
 {}\times
 F_{cbe}(s_1,\chi)
 \overline{F_{cbe}(\overline{s_2},\chi)}
 \widetilde{\mathcal W}_2
 \left(s_1,s_2;\frac{X}{cbgQ}\right)
 \dif s_1\dif s_2,\nonumber
\end{align}
where $\beta>1$; the character variable has
been relabelled by $\chi\mapsto\overline\chi$. Define
\[
 F_M(s,\chi)=\sum_{n=1}^\infty
 \frac{\lambda_{\pi}(Mn)\overline{\chi(n)}}{n^s}.
\]
The Euler product gives
\begin{align*}
 F_M(s,\chi)
 =L(s,\pi\otimes\overline\chi)E_M(s,\chi).
\end{align*}
If $p^\ell\Vert M$, its local factor is
\begin{align}\label{eq:EA-local-factor}
 E_{p^\ell}(s,\chi)
 =\frac{\displaystyle
  \sum_{k\geq0}\lambda_{\pi}(p^{\ell+k})\overline{\chi(p)}^{k}p^{-ks}}
 {\displaystyle
  \sum_{k\geq0}\lambda_{\pi}(p^k)\overline{\chi(p)}^{k}p^{-ks}}.
\end{align}
For $p$ dividing the modulus of $\chi$, the convention $\chi(p)=0$
reduces \eqref{eq:EA-local-factor} to $\lambda_{\pi}(p^\ell)$.  Thus the same formula
covers these bad primes.

\begin{lem}\label{lem:EA-average}
Uniformly for real $t$, every Dirichlet character $\chi$, and $Y\geq2$,
\begin{align}\label{eq:EA-average}
 \sum_{M\leq Y}\frac{|E_M(1/2+it,\chi)|^2}{M}
 \ll_{\pi,\varepsilon}(Y(2+|t|))^\varepsilon.
\end{align}
The estimate remains valid if $M$ is restricted by finitely many divisor or
coprimality conditions.
\end{lem}

\begin{proof}
Write
\[
 L_p(s,\pi)^{-1}
 =\sum_{u=0}^{m}\mu_\pi(p^u)p^{-us},
\]
inserting zero coefficients when the inverse local factor has degree less
than \(m\).  If \(p^\ell\Vert M\), put
\[
 Z=\overline{\chi(p)}p^{-1/2-it},
 \qquad |Z|\leq p^{-1/2}.
\]
Then \eqref{eq:EA-local-factor} gives
\[
 E_{p^\ell}(1/2+it,\chi)
 =
 \left(\sum_{u=0}^{m}\mu_\pi(p^u)Z^u\right)
 \left(\sum_{k\geq0}\lambda_{\pi}(p^{\ell+k})Z^k\right).
\]
Since
\[
 \sum_{u=0}^{m}
 \mu_\pi(p^u)\lambda_{\pi}(p^{k-u})=0
 \qquad(k\geq1),
\]
with \(\lambda_{\pi}(p^v)=0\) for \(v<0\), the coefficients of \(Z^v\) vanish for
\(v\geq m\).  Hence
\[
 E_{p^\ell}(1/2+it,\chi)
 =
 \sum_{v=0}^{m-1}
 \left(
  \sum_{u=0}^{v}
  \mu_\pi(p^u)\lambda_{\pi}(p^{\ell+v-u})
 \right)Z^v.
\]
Cauchy--Schwarz therefore gives
\begin{align}\label{eq:EA-local-square-bound}
 |E_{p^\ell}(1/2+it,\chi)|^2
 \ll_\pi
 \sum_{v=0}^{m-1}\sum_{u=0}^{v}
 \frac{|\mu_\pi(p^u)|^2
       |\lambda_{\pi}(p^{\ell+v-u})|^2}{p^v}.
\end{align}

Let \(1<w\le2\) and set
\[
 \sigma=1+\frac{w-1}{2m}.
\]
For \(\ell\geq1\) and \(0\leq u\leq v\leq m-1\),
\[
 \ell w+v-\bigl((\ell+v-u)+u\bigr)\sigma
 =(w-1)\left(\ell-\frac{\ell+v}{2m}\right)
 \geq\frac{w-1}{2}.
\]
Consequently,
\[
 \frac{1}{p^{\ell w+v}}
 \leq
 \frac{1}{p^{(\ell+v-u)\sigma}p^{u\sigma}}.
\]
Summing \eqref{eq:EA-local-square-bound} over \(\ell\), and noting that
each pair \((u,\ell+v-u)\) occurs only \(O_m(1)\) times, shows that
there is a fixed integer \(C=C(\pi)\) such that
\begin{align*}
 1+\sum_{\ell\geq1}
 \frac{|E_{p^\ell}(1/2+it,\chi)|^2}{p^{\ell w}}
 \leq
 \left(1+\sum_{k\geq1}
       \frac{|\lambda_{\pi}(p^k)|^2}{p^{k\sigma}}\right)^C
 \left(1+\sum_{u=1}^{m}
       \frac{|\mu_\pi(p^u)|^2}{p^{u\sigma}}\right)^C.
\end{align*}

Lemma~\ref{corollary:Jiang} gives
\[
 \sum_{n\geq1}\frac{|\lambda_\pi(n)|^2}{n^\sigma}
 \ll_\pi(\sigma-1)^{-C},
\]
after enlarging \(C\) if necessary.  To treat the reciprocal coefficients,
use
\[
 \mu_\pi=\mu_{\pi_1}\ast\cdots\ast\mu_{\pi_r}.
\]
Cauchy--Schwarz and \eqref{ineq:Coefficient&Reciprocal} imply, again after
enlarging \(C\),
\[
 1+\sum_{u=1}^{m}
 \frac{|\mu_\pi(p^u)|^2}{p^{u\sigma}}
 \leq
 \prod_{j=1}^r
 L_p(\sigma,\pi_j\times\widetilde\pi_j)^C.
\]
It follows that
\[
 \sum_{n\geq1}\frac{|\mu_\pi(n)|^2}{n^\sigma}
 \ll_\pi
 \prod_{j=1}^r
 L(\sigma,\pi_j\times\widetilde\pi_j)^C
 \ll_\pi(\sigma-1)^{-C}.
\]

Since \(E_M(s,\chi)\) is multiplicative in \(M\), the preceding estimates
yield
\begin{align*}
 \sum_{M\geq1}
 \frac{|E_M(1/2+it,\chi)|^2}{M^w}
 &=
 \prod_p\left(
  1+\sum_{\ell\geq1}
  \frac{|E_{p^\ell}(1/2+it,\chi)|^2}{p^{\ell w}}
 \right)\\
 &\ll_\pi
 \left(
  \sum_{n\geq1}\frac{|\lambda_\pi(n)|^2}{n^\sigma}
 \right)^C
 \left(
  \sum_{n\geq1}\frac{|\mu_\pi(n)|^2}{n^\sigma}
 \right)^C\\
 &\ll_\pi(w-1)^{-2C}.
\end{align*}
The argument may initially be applied to finite Euler products and then
completed by monotone convergence.

Taking
\[
 w=1+\frac{1}{\log(2Y)}
\]
and applying Rankin's trick, we obtain
\begin{align*}
 \sum_{M\leq Y}
 \frac{|E_M(1/2+it,\chi)|^2}{M}
 &\leq
 Y^{w-1}
 \sum_{M\geq1}
 \frac{|E_M(1/2+it,\chi)|^2}{M^w}\\
 &\ll_\pi
 \bigl(\log(2Y)\bigr)^{2C}
 \ll_{\pi,\varepsilon}Y^\varepsilon.
\end{align*}
The bound is uniform in \(\chi\) and \(t\), and it is stronger than
\eqref{eq:EA-average}.  Any divisor or coprimality restriction on \(M\)
only removes nonnegative terms from the sum.
\end{proof}

Since $F_M(s,\chi)$ is entire, we may shift the two contours in
\eqref{eq:G2-Mellin} to $\operatorname{Re}(s_1)=\operatorname{Re}(s_2)=1/2$.
Put
\[
 h=cg,\qquad M=cb e.
\]
The support of \(\Psi\) implies \(M\ll X\), while the support statement in
Proposition~\ref{MultiMellinBound2} implies
\[
 hb\ll \frac{X}{Q}.
\]
In particular, \(c\mid h\) and \(e\ll X/(cb)\).

Write
\[
 \mathcal K_{h,b}(t_1,t_2)
 =
 \left|
 \widetilde{\mathcal W}_2
 \left(
  \tfrac12+it_1,\tfrac12+it_2;
  \frac{X}{hb Q}
 \right)
 \right|.
\]
Write $L_\chi(t)=L(1/2+it,\pi\otimes\overline\chi)$ for $\chi\bmod h$.
Applying $2|uv|\le|u|^2+|v|^2$ to the two $F_M$ factors in
\eqref{eq:G2-Mellin} and using $\phi(h)\gg_\varepsilon h^{1-\varepsilon}$
gives terms of the form
\[
 \frac{X^{1+\varepsilon}}{hcb e}
 |L_\chi(t)|^2
 |E_{cb e}(1/2+it,\chi)|^2
 \mathcal K_{h,b}(t_1,t_2).
\]
Here $t=t_1$ or $t=-t_2$. In the second term we may change
$(t_1,t_2)$ to $(-t_1,-t_2)$: the ranges in both variables are
symmetric and the real weight gives
$\mathcal K_{h,b}(-t_1,-t_2)=\mathcal K_{h,b}(t_1,t_2)$.

 For a dyadic number \(B\), put
\[
 \mathcal{B}_{h,B}(t,\chi)
 =
 \sum_{b\asymp B}
 \sum_{c\mid h}
 \sum_{\substack{e\ll X/(cb)\\(e,h)=1}}
 \frac{
  |E_{cb e}(1/2+it,\chi)|^2
 }{cb e}.
\]
Writing \(M=cb e\), we have
\[
 \mathcal B_{h,B}(t,\chi)
 =
 \sum_{M\ll X}
 \frac{|E_M(1/2+it,\chi)|^2}{M}
 R_{h,B}(M),
\]
where
\[
 R_{h,B}(M)
 =
 \#\left\{
 (c,b,e):
 c\mid h,\ b\asymp B,\ M=cb e,\ (e,h)=1
 \right\}.
\]
For each \(M\),
\[
 R_{h,B}(M)\leq \tau(h)\tau(M)
 \ll_\varepsilon (hM)^\varepsilon.
\]
Consequently, Lemma~\ref{lem:EA-average} gives
\begin{align*}
 \mathcal B_{h,B}(t,\chi)
 &\ll_{\pi,\varepsilon}
 (hX)^\varepsilon
 \sum_{M\ll X}
 \frac{|E_M(1/2+it,\chi)|^2}{M} \notag\\
 &\ll_{\pi,\varepsilon}(hX)^\varepsilon,
\end{align*}

For \(S>1\), set
\[
 \mathcal E(S)=[S,2S]\cup[-2S,-S],
 \qquad
 \mathcal E(1)=[-2,2].
\]
Decompose
\[
 h\asymp H,\qquad
 b\asymp B,\qquad
 t_k\in\mathcal E(S_k),\qquad
 1\leq S_1\leq S_2,
\]
where \(H,B,S_1,S_2\) run over dyadic numbers and
\[
 HB\ll \frac{X}{Q}.
\]
The ranges with \(S_2<S_1\) are treated by interchanging \(t_1\) and
\(t_2\).

For fixed \(h\asymp H\), define
\[
 \mathcal K^*_{h,B}(t_1,t_2)
 =
 \sup_{b\asymp B}
 \mathcal K_{h,b}(t_1,t_2).
\]
Using the uniform estimate proved in
Lemma~\ref{lem:EA-average}, and absorbing the resulting
\(X^\varepsilon\)-factor by changing \(\varepsilon\), we obtain
\begin{align}\label{eq:G2-reduced}
 |\mathscr B_{2,\pi}(Q,X)|
 \ll_{\pi,\varepsilon}
 X^{1+\varepsilon}
 \sum_{\substack{H,B\ {\rm dyadic}\\HB\ll X/Q}}
 \ \sum_{\substack{S_1,S_2\ {\rm dyadic}\\1\leq S_1\leq S_2}}
 \mathcal J(H,B;S_1,S_2),
\end{align}
where
\begin{align*}
 \mathcal J(H,B;S_1,S_2)
 \ll_{\pi,\varepsilon}
 H^\varepsilon
 \sum_{h\asymp H}\frac1h
 \sum_{\substack{\chi\bmod h\\\chi\ne\chi_0}}
 \int_{\mathcal E(S_1)}
 \int_{\mathcal E(S_2)}
 \Bigl(
  |L_\chi(t_1)|^2+|L_\chi(t_2)|^2
 \Bigr)
 \mathcal K^*_{h,B}(t_1,t_2)
 \dif t_2\dif t_1.
\end{align*}

Put
\[
 U=\frac{X}{HBQ}.
\]
For \(h\asymp H\), \(b\asymp B\) and every integer \(P\geq1\),
Proposition~\ref{MultiMellinBound2} gives uniform bound
\begin{align}\label{eq:K-dyadic-bound}
 \mathcal K^*_{h,B}(t_1,t_2)
 \ll_P
 \frac{U^{P-1}}{S_2^P}
 \qquad
 \left(
  t_k\in\mathcal E(S_k),\ S_1\leq S_2
 \right).
\end{align}

We next apply Lemma~\ref{lem:hybridlargesieve} separately to the two
\(L\)-function terms.  For the term containing \(|L_\chi(t_1)|^2\), the
\(t_2\)-integration contributes \(O(S_2)\), whereas for the term containing
\(|L_\chi(t_2)|^2\), the \(t_1\)-integration contributes \(O(S_1)\).
Together with \eqref{eq:K-dyadic-bound}, this gives
\begin{align*}
 \mathcal J(H,B;S_1,S_2)
 \ll_{\pi,P,\varepsilon}
 \frac{U^{P-1}}{HS_2^P}
 (HS_2)^\varepsilon
 \Bigl\{&
 S_2\bigl(H^2S_1+(HS_1)^{m/2}\bigr)\\
 &+
 S_1\bigl(H^2S_2+(HS_2)^{m/2}\bigr)
 \Bigr\}.
\end{align*}
If \(m\geq2\), then \(S_1\leq S_2\) implies
\[
 S_2(HS_1)^{m/2}
 \leq
 S_1(HS_2)^{m/2}.
\]
If \(m=1\), the term
\(S_2(HS_1)^{1/2}\) is bounded by \(H^2S_1S_2\), since
\(H,S_1\geq1\).  Thus, in every case,
\begin{align}\label{eq:G2-dyadic}
 \mathcal J(H,B;S_1,S_2)
 \ll_{\pi,P,\varepsilon}
 \frac1H
 \left(
  H^2S_2+(HS_2)^{m/2}
 \right)
 (HS_2)^\varepsilon
 S_1
 \frac{U^{P-1}}{S_2^P}.
\end{align}

Set
\[
 \gamma=\frac{m}{2}.
\]
For \(S_2\leq U\), use \eqref{eq:G2-dyadic} with \(P=1\).
Since the dyadic sum of \(S_1\) over \(S_1\leq S_2\) is \(O(S_2)\),
the contribution from this range is
\[
 \ll_{\pi,\varepsilon}
 \frac{H^2U+(HU)^\gamma}{H}
 (HU)^\varepsilon.
\]
For \(S_2>U\), choose the integer
\[
 P=\lceil\gamma\rceil+2.
\]
Then \(P>2\) and \(P>\gamma+1\), so both dyadic \(S_2\)-sums converge
geometrically.  Their contribution is again
\[
 \ll_{\pi,\varepsilon}
 \frac{H^2U+(HU)^\gamma}{H}
 (HU)^\varepsilon.
\]
It follows that
\begin{align}\label{eq:G2-time-summed}
 \sum_{\substack{S_1,S_2\ {\rm dyadic}\\1\leq S_1\leq S_2}}
 \mathcal J(H,B;S_1,S_2)
 \ll_{\pi,\varepsilon}
 \frac{H^2U+(HU)^{m/2}}{H}
 (HU)^\varepsilon.
\end{align}

Since \(U=X/(HBQ)\), the right-hand side of
\eqref{eq:G2-time-summed} is
\[
 \ll_{\pi,\varepsilon}
 \left\{
  \frac{X}{BQ}
  +
  \frac{(X/Q)^{m/2}}{H B^{m/2}}
 \right\}
 (HU)^\varepsilon.
\]
The dyadic \(B\)-sum in the first term is convergent, while the corresponding
\(H\)-sum contributes only \(O(\log X)\).  The dyadic \(H\)- and \(B\)-sums
in the second term are both convergent.  Absorbing the logarithmic factor
into \(X^\varepsilon\), \eqref{eq:G2-reduced} therefore gives
\begin{align*}
 \mathscr B_{2,\pi}(Q,X)
 \ll_{\pi,\Phi,\Psi,\varepsilon}
 X^{1+\varepsilon}
 \left\{
  \frac XQ+
  \left(\frac XQ\right)^{m/2}
 \right\}.
\end{align*}
This proves the unconditional assertion of
Proposition~\ref{finalg2}.

Under GRH, Lemma~\ref{lem:hybridlargesieve} contains only the term
\(H^2S\).  Repeating the preceding argument gives
\[
 \sum_{S_1,S_2}
 \mathcal J(H,B;S_1,S_2)
 \ll_{\pi,\varepsilon}
 HU(HU)^\varepsilon
 =
 \frac{X}{BQ}
 \left(\frac{X}{BQ}\right)^\varepsilon.
\]
Summing over the dyadic \(H,B\), with the logarithmic \(H\)-sum absorbed
into \(X^\varepsilon\), yields
\begin{align*}
 \mathscr B_{2,\pi}(Q,X)
 \ll_{\pi,\Phi,\Psi,\varepsilon}
 \frac{X^{2+\varepsilon}}{Q}.
\end{align*}
Comparison with \(QX\) completes the proof of
Proposition~\ref{finalg2}.

\section{Proof of Lemma~\ref{analyticcontinuationofTk}}
\label{sec:proof-analyticcontinuationofTk}

We begin in the region
$\mathrm{Re}(s_1),\mathrm{Re}(s_2)>1$ and $\mathrm{Re}(z)<0$, where all rearrangements are
absolutely convergent.

The local calculations below are made at the unramified primes
\(p\nmid q_\pi\). At \(p\mid q_\pi\) we retain the actual local factors
in the quotient \(J_p\) defined below. Write \(w=s_1+s_2\), and take
\(\vartheta=\theta\) unconditionally and \(\vartheta=0\) under GRC.
In the corresponding region of Lemma~\ref{analyticcontinuationofTk},
\(\Re s_i-\vartheta\geq 1/3+\varepsilon\) and
\(\min\{\Re w,\Re(w-z)\}-2\vartheta\geq 2/3\).
The local parameter bounds of Section~2 and the geometric sums in
\eqref{eq:Tp-exact} thus give \(|J_p|\leq K_{m,\varepsilon}\) uniformly
in \(p\) and in the three variables. Hence the ramified part satisfies
\[
 \prod_{p\mid q_\pi}|J_p|
 \leq K_{m,\varepsilon}^{\omega(q_\pi)}
 \ll_{m,\varepsilon} q_\pi^\varepsilon
 \leq \mathcal C(\pi)^\varepsilon,
\]
where \(\omega(q_\pi)\) counts the distinct prime divisors of \(q_\pi\).

Except for symbols with an established meaning, the local factors introduced
in this section are denoted in alphabetical order.  Each defining display is
tagged by the corresponding local-factor symbol for ease of reference.

\subsection{The auxiliary divisor sum}

For fixed $d,N_1,N_2$, put $t=cg$.  Since $c$ is squarefree,
\begin{align}\label{eq:ar-sum-initial}
 \sum_{g\geq1}\sum_{\substack{c\mid d\\(N_1N_2,cg)=1}}
 \frac{\mu(c)g^z}{\phi(cg)}
 =\sum_{\substack{c\mid d\\(c,N_1N_2)=1}}
 \frac{\mu(c)}{c^z}
 \sum_{\substack{t\equiv0\pmod c\\(t,N_1N_2)=1}}
 \frac{t^z}{\phi(t)}.
\end{align}
Define
\begin{align}
 A_p(z)
 &=1+\sum_{\rho\geq1}\frac{p^{\rho z}}{\phi(p^\rho)}
   =\frac{1+p^{z-1}/(p-1)}{1-p^{z-1}},\tag{$A_p$}\\
 B_p(z)
 &=1-p^{-z}\bigl(1-A_p(z)^{-1}\bigr)
   =1-\frac1{p-1+p^{z-1}}.\tag{$B_p$}
\end{align}
The inner sum in \eqref{eq:ar-sum-initial} is multiplicative and equals
\begin{align*}
 \prod_pA_p(z)
 \prod_{p\mid c}\bigl(1-A_p(z)^{-1}\bigr)
 \prod_{p\mid N_1N_2}A_p(z)^{-1}.
\end{align*}
Summing over $c$ gives the exact identity
\begin{align}\label{eq:ar-sum-factorization}
 \sum_{g\geq1}\sum_{\substack{c\mid d\\(N_1N_2,cg)=1}}
 \frac{\mu(c)g^z}{\phi(cg)}
 =\prod_pA_p(z)
  \prod_{p\mid N_1N_2}A_p(z)^{-1}
  \prod_{\substack{p\mid d\\p\nmid N_1N_2}}B_p(z).
\end{align}
Moreover,
\begin{align}\label{eq:zeta-factor}
 \prod_pA_p(z)
 =\zeta(1-z)\mathcal P_0(z),\hspace{2em}
 \mathcal P_0(z)=\prod_p\left(1+\frac{p^{z-1}}{p-1}\right),
\end{align}
and $\mathcal P_0$ converges absolutely for $\mathrm{Re}(z)<1$.

\subsection{The coefficient sums}

At a fixed prime $p$ write
\[
 e=v_p(d),\hspace{2em}u=v_p(N_1),\hspace{2em}v=v_p(N_2).
\]
The condition $(N_1,N_2)=1$ is $\min(u,v)=0$.  If
$\xi=v_p(c)\in\{0,1\}$ and $\rho=v_p(g)$, the further conditions are
$\xi\leq e$ and
\[
 u+v>0\quad\Longrightarrow\quad \xi=\rho=0.
\]
Consequently the local Euler factor of $T$ is the exact sum
\begin{align}\label{eq:Tp-exact}
 T_p(s_1,s_2;z)
 =\sum_{\substack{e,u,v\geq0\\\min(u,v)=0}}
 \sum_{\substack{\xi\in\{0,1\},\ \rho\geq0\\
                  \xi\leq e\\
                  (\xi+\rho)(u+v)=0}}
 \frac{\overline{\lambda_{\pi}(p^{e+u})}\lambda_{\pi}(p^{e+v})
       \mu(p^\xi)p^{(e+\rho)z}}
      {p^{(e+u)s_1+(e+v)s_2}\phi(p^{\xi+\rho})}.
\end{align}

The decomposition of \eqref{eq:Tp-exact} is most transparent after the
following definitions.  For $e\geq0$ introduce the coefficient and tail
factors
\begin{align}
 C_{p,e}(s;\pi)
 &=\frac{\lambda_{\pi}(p^e)}{p^{es}},\tag{$C_{p,e}$}\\
 D_{p,e}(s;\pi)
 &=L_p(s,\pi)-\sum_{k=0}^{e}\frac{\lambda_{\pi}(p^k)}{p^{ks}}
   =\sum_{k=e+1}^{\infty}\frac{\lambda_{\pi}(p^k)}{p^{ks}}.\tag{$D_{p,e}$}
\end{align}
The following three factors separate the cases $e=0$, $u=v=0$, and
$e(u+v)>0$:
\begin{align}
 E_p(s_1,s_2;z)
 &=A_p(z)^{-1}\left\{
   \overline{D_{p,0}(\overline{s_1};\pi)}
   +D_{p,0}(s_2;\pi)\right\},\tag{$E_p$}\\[1.2ex]
 F_{p,e}(s_1,s_2;z)
 &=B_p(z)\overline{C_{p,e}(\overline{s_1};\pi)}
              C_{p,e}(s_2;\pi),\qquad e\geq1,\tag{$F_{p,e}$}\\[1.2ex]
 G_{p,e}(s_1,s_2;z)
 &=A_p(z)^{-1}\left\{
 C_{p,e}(s_2;\pi)\overline{D_{p,e}(\overline{s_1};\pi)}
 +\overline{C_{p,e}(\overline{s_1};\pi)}D_{p,e}(s_2;\pi)
 \right\},\qquad e\geq1.\tag{$G_{p,e}$}
\end{align}
Here $E_p$ contains the terms in which $p$ divides exactly one of $N_1,N_2$ but
not $d$; $F_{p,e}$ contains $p^e\Vert d$ with $p\nmid N_1N_2$; and $G_{p,e}$
contains $p^e\Vert d$ together with a positive exponent in one of $N_1,N_2$.
Using \eqref{eq:ar-sum-factorization} and $\min(u,v)=0$, the normalized local
factor is therefore
\begin{align}\label{eq:Tp-truncated-series}
 A_p(z)^{-1}T_p(s_1,s_2;z)
 =1+E_p(s_1,s_2;z)
  +\sum_{e\geq1}p^{ez}
   \left\{F_{p,e}(s_1,s_2;z)+G_{p,e}(s_1,s_2;z)\right\}.
\end{align}
The conjugated notation in the first variable is used only to reproduce the
coefficients $\overline{\lambda_{\pi}(p^k)}$; explicitly,
\[
 \overline{D_{p,e}(\overline{s_1};\pi)}
 =\sum_{k=e+1}^{\infty}\frac{\overline{\lambda_{\pi}(p^k)}}{p^{ks_1}}.
\]
Thus the shift by $e=v_p(d)$ is incorporated in the lower endpoint
$k=e+1$ of the tail.

The Satake bound gives, uniformly for $\sigma>\theta$,
\begin{align}\label{eq:local-tail-bound}
 C_{p,e}(\sigma+it;\pi)&\ll_{\pi}(e+1)^{m-1}
        p^{-e(\sigma-\theta)},\\
 D_{p,e}(\sigma+it;\pi)&\ll_{\pi,\sigma}
        (e+2)^{m-1}p^{-(e+1)(\sigma-\theta)}.\nonumber
\end{align}
Also, locally uniformly for $\mathrm{Re}(z)<1$,
\begin{align*}
 A_p(z)^{-1}=1+O\left(p^{\mathrm{Re}(z)-1}\right),
 \qquad
 B_p(z)=1+O(p^{-1}).
\end{align*}
Write \(\sigma_k=\Re(s_k)\), \(\sigma=\min(\sigma_1,\sigma_2)\),
\(\tau=\Re(z)\) and \(b=\sigma_1+\sigma_2-2\theta-\tau\).
In the region of Lemma~\ref{analyticcontinuationofTk}, \(b\ge 2/3\).
It follows that
\begin{align*}
 E_p(s_1,s_2;z)
 ={}&\frac{\overline{\lambda_{\pi}(p)}}{p^{s_1}}
    +\frac{\lambda_{\pi}(p)}{p^{s_2}}
    +O_\pi\left(p^{-2(\sigma-\theta)}
      +p^{\mathrm{Re}(z)-1-(\sigma-\theta)}\right),\\
 p^zF_{p,1}(s_1,s_2;z)
 ={}&B_p(z)\frac{|\lambda_{\pi}(p)|^2}{p^{s_1+s_2-z}},\nonumber\\
 \sum_{e\geq1}\bigl|p^{ez}G_{p,e}(s_1,s_2;z)\bigr|
 \ll{}&p^{3\theta-\mathrm{Re}(s_1+s_2)-\sigma+\mathrm{Re}(z)}.\nonumber
\end{align*}
Indeed, \eqref{eq:local-tail-bound} gives
\[
 |p^{ez}G_{p,e}|\ll_\pi (e+2)^{2m-2}p^{-(\sigma-\theta)-eb}.
\]
Summing over \(e\ge1\) costs \(O_m(p^{-b})\) after factoring out
\(p^{-(\sigma-\theta)}\).

\subsection{Extraction of the automorphic factors}

Define the remaining local factors by
\begin{align}
 H_p(s_1,s_2)
 &=L_p(s_1,\widetilde\pi)L_p(s_2,\pi),\tag{$H_p$}\\
 I_p(s_1,s_2;z)
 &=\frac{\displaystyle\prod_{h,j}
     L_p(s_1+s_2-z,\pi_h\times\widetilde\pi_j)}
    {\displaystyle\prod_{h,j}
     L_p(s_1+s_2,\pi_h\times\widetilde\pi_j)},\tag{$I_p$}\\
 J_p(s_1,s_2;z)
 &=\frac{A_p(z)^{-1}T_p(s_1,s_2;z)}
        {H_p(s_1,s_2)I_p(s_1,s_2;z)}.\tag{$J_p$}
\end{align}
The first two have the expansions, with $w=s_1+s_2$,
\begin{align}\label{eq:HI-expansions}
 H_p(s_1,s_2)
 ={}&1+\frac{\overline{\lambda_{\pi}(p)}}{p^{s_1}}
      +\frac{\lambda_{\pi}(p)}{p^{s_2}}
      +\frac{|\lambda_{\pi}(p)|^2}{p^w}+\mathcal H_p,\\
 I_p(s_1,s_2;z)
 ={}&1+|\lambda_{\pi}(p)|^2
 \left(\frac1{p^{w-z}}-\frac1{p^w}\right)+\mathcal I_p.\nonumber
\end{align}
Here $\mathcal H_p$ and $\mathcal I_p$ consist of terms of local degree at
least two.  The terms $|\lambda_{\pi}(p)|^2p^{-w}$ in \eqref{eq:HI-expansions} cancel,
so the first-order local terms of $H_pI_p$ agree with those of
\eqref{eq:Tp-truncated-series}, in view of
$B_p(z)=1+O(p^{-1})$.  It remains to prove the absolute convergence of
$\prod_pJ_p$.

We give the convergence argument in local form.  Write
\[
 \sigma_k=\mathrm{Re}(s_k),\qquad \sigma=\min\{\sigma_1,\sigma_2\},
 \qquad w=s_1+s_2.
\]
The one-sided terms represented by $E_p$ are collected into the full factor
$H_p$.  The terms $F_{p,e}$ with
$u=v=0$ form a square-coefficient series at $w-z$, whereas the omission of
$u,v>0$ simultaneously produces the reciprocal square term at $w$.
The local forms of \eqref{eq:H-factor} and \eqref{eq:H1-factor} give
\begin{align}\label{eq:local-square-extraction}
 1+\sum_{e\geq1}\frac{|\lambda_{\pi}(p^e)|^2}{p^{e(w-z)}}
 &=H_{\pi,p}(w-z)
   \prod_{h,j}L_p(w-z,\pi_h\times\widetilde\pi_j),\\
 1-\frac{|\lambda_{\pi}(p)|^2}{p^w}
 &=H_{\pi,p}^{\sharp}(w)
   \prod_{h,j}L_p(w,\pi_h\times\widetilde\pi_j)^{-1}.
 \nonumber
\end{align}
The Euler products of the $H$-factors in
\eqref{eq:local-square-extraction} converge throughout the region considered
below and are absorbed into $\mathcal P$.

After the factors $H_p$ and $I_p$ have been removed, the local quotient
$J_p$ has three kinds of remaining terms.  First, the factors $G_{p,e}$ are
bounded, using \eqref{eq:local-tail-bound}, by
\begin{align}\label{eq:mixed-local-error}
 O_\pi\left(
 p^{3\theta-\mathrm{Re}(w)-\sigma+\mathrm{Re}(z)}
 \right).
\end{align}
Second, replacing $A_p(z)^{-1}$ by $1$ in the one-sided factors contributes
\begin{align}\label{eq:A-local-error}
 O_\pi\left(
 p^{\mathrm{Re}(z)-1-\sigma+\theta}
 \right).
\end{align}
Third, replacing $B_p(z)$ by $1$ in the square-coefficient series contributes
\begin{align}\label{eq:B-local-error}
 O_\pi\left(
 p^{-1-\mathrm{Re}(w-z)+2\theta}
 \right).
\end{align}
All further terms contain an extra factor
$p^{-(\sigma-\theta)}$ or a prime power coefficient of order at least two.
The former sum geometrically by \eqref{eq:local-tail-bound}; the latter are
summed by Theorem~\ref{Jiang} and
Lemma~\ref{corollary:Jiang}.

Suppose now that
\[
 \sigma_1,\sigma_2\geq\frac13+\theta+\varepsilon,
 \qquad \mathrm{Re}(z)\leq2\varepsilon.
\]
Then the exponent in \eqref{eq:mixed-local-error} is at most
$-1-\varepsilon$, the exponent in \eqref{eq:A-local-error} is at most
$-4/3+\varepsilon$, and the exponent in \eqref{eq:B-local-error} is at most
$-5/3$.  On decreasing $\varepsilon$ in the initial statement if necessary,
there exists $\eta_0>0$ such that the normalized local quotient is
\[
 1+O_{\pi,\varepsilon}(p^{-1-\eta_0}).
\]
The estimate is locally uniform in all three variables.  Therefore
\begin{align*}
 \mathcal P(s_1,s_2,z)
 =\mathcal P_0(z)\prod_pJ_p(s_1,s_2;z)
\end{align*}
is a bounded, absolutely convergent Euler product.

Combining \eqref{eq:zeta-factor} with the extracted local factors proves
\[
 T(s_1,s_2;z)
 =\zeta(1-z)
 \frac{\displaystyle
  \prod_hL(s_1,\widetilde\pi_h)
  \prod_jL(s_2,\pi_j)
  \prod_{h,j}L(s_1+s_2-z,\pi_h\times\widetilde\pi_j)}
 {\displaystyle\prod_{h,j}L(s_1+s_2,\pi_h\times\widetilde\pi_j)}
 \mathcal P(s_1,s_2,z).
\]
This gives the asserted continuation.  Under GRC the local parameter bounds hold with $\theta=0$,
which enlarges the region to
$\mathrm{Re}(s_1),\mathrm{Re}(s_2)\geq1/3+\varepsilon$ and
$\mathrm{Re}(z)\leq2\varepsilon$.

Specializing the coefficients in the definition of \(T\) to \(d_k\)
gives the series \(T_k\) of Rodgers and Soundararajan \cite{RS}.
The two standard \(L\)-factors then become \(\zeta(s_1)^k\) and
\(\zeta(s_2)^k\), whose poles give rise to the residues in their
contour shift.

\end{document}